\documentclass[10pt,twocolumn,letterpaper]{article}

\usepackage{cvpr}
\usepackage{times}
\usepackage{epsfig}
\usepackage{graphicx}
\usepackage{amsmath}
\usepackage{amssymb}
\usepackage[acronym]{glossaries}

\usepackage{caption}
\usepackage{subcaption}
\usepackage{algorithm}
\usepackage{algpseudocode}
\usepackage{nomencl}
\makenomenclature
\usepackage{etoolbox}
\renewcommand\nomgroup[1]{%
  \item[\bfseries
  \ifstrequal{#1}{V}{Variables}{}
  \ifstrequal{#1}{A}{Acronyms}{}%
]}

\usepackage[official]{eurosym}

\usepackage{natbib}

\newcommand{\bx}{x}
\newcommand{\bxh}{\hat{x}}
\newcommand{\bxt}{\tilde{x}}
\newcommand{\bA}{SA}
\newcommand{\bEA}{EA}

\usepackage[breaklinks=true,bookmarks=false]{hyperref}

\cvprfinalcopy 

\def\cvprPaperID{****} 

\newacronym{mrp}{MRP}{Material Requirement Planning}
\newacronym{otif}{OTIF}{On-time in-full}

\begin{document}

\title{Stochastic Optimization of Inventory at Large-scale Supply Chains}

\author{Zhaoyang Larry Jin, \quad Mehdi Maasoumy, \quad Yimin Liu, \quad Zeshi Zheng, \quad Zizhuo Ren \\
C3 AI\\
{\tt\small \{larry.jin, mehdi.maasoumy, yimin.liu, zeshi.zheng, zizhuo.ren\}@c3.ai}}

\maketitle

\begin{abstract}
Today’s global supply chains face growing challenges due to rapidly changing market conditions, increased network complexity and inter-dependency, and dynamic uncertainties in supply, demand, and other factors. To combat these challenges, organizations employ Material Requirements Planning (MRP) software solutions to set inventory stock buffers – for raw materials, work-in-process goods, and finished products – to help them meet customer service levels. However, holding excess inventory further complicates operations and can lock up  millions of dollars of capital that could be otherwise deployed. Furthermore, most commercially available MRP solutions fall short in considering uncertainties and do not result in optimal solution for modern enterprises.

At C3 AI, we fundamentally reformulate the inventory management problem as a constrained stochastic optimization. We then propose a simulation-optimization framework that minimizes inventory and related costs while maintaining desired service levels. The framework’s goal is to find the optimal reorder parameters  that minimize costs subject to a pre-defined service-level constraint and all other real-world operational constraints. These optimal reorder parameters can be fed back into an MRP system to drive optimal order placement, or used to place optimal orders directly. This approach has proven successful in reducing inventory levels by 10-35 percent, resulting in hundreds of millions of dollars of economic benefit for major enterprises at a global scale.

\end{abstract}

\section{Introduction}
\label{section:introduction}

Global supply chain faced enormous challenges due to uncertainties introduced by the shocks in both supply side and demand side during the COVID-19 pandemic. A supply shock hit as the pandemic first emerged from mainland China and disrupted global supply. Over the ensuing months, global shutdowns and stay-at-home orders triggered a second shock, this time to the global demand. These two waves of unprecedented disruptions to the global supply chains exposed the inherent volatility of such networks. With this volatility to be continued in the post-pandemic era, a solution is required to handle uncertainties in the global supply chains. While battling the uncertainties by increasing capital spent on raising the inventory level, the enterprises begin to pay more attention at the concept of lean manufacturing which includes maximizing the resilience of the inventory availability while minimizing the amount of capital expenditure on inventory held.

\begin{figure}[t]
  \centering
  \includegraphics[width=\linewidth]{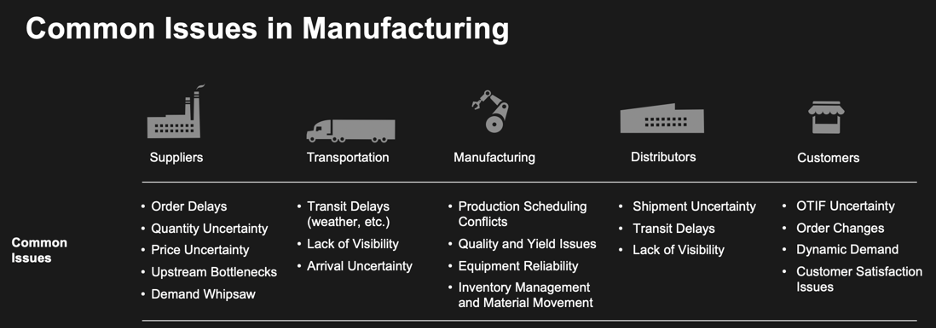}
  \caption{Common issues in manufacturing ranging from suppliers to customers}
  \label{fig:common_issues}
\end{figure}

Current practice in today’s supply chains suffer from numerous issues covering aspects including suppliers, transportation, manufacturing, distribution and customers (Fig.~\ref{fig:common_issues}). Perturbations in each of the aspects are hardly specific only to itself. The network nature of the supply chains makes any issue in each of these aspects easily propagate and amplify upstream or downstream. On one direction, dynamic demand and order changes would impact production schedules and inventory management at manufacturing sites or distribution centers. In addition, available inventory and production requirements would inform supplier orders. The equipment reliability, product quality issues, and shipment uncertainties would impact part orders. On the other direction, supplier uncertainties and delays would affect production schedules, while order delays and quality issue would risk customer satisfaction. Furthermore, material movement delays would strain the entire supply chain. The Lack of end-to-end operational visibility would lead to delays across supply chains, excess inventory levels, quality issues, poor OTIF, or low customer satisfaction.

\begin{figure}[!htp]
  \centering
  \includegraphics[width=\linewidth]{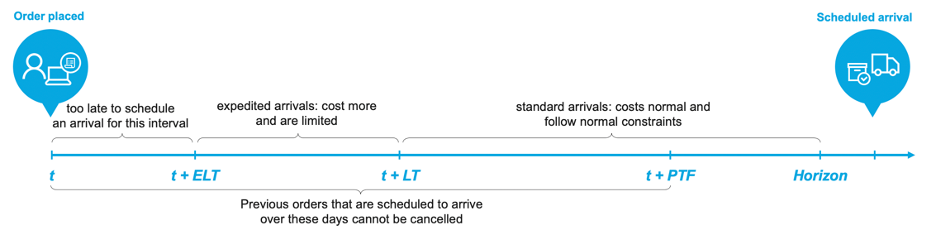}
  \caption{MRP follows a deterministic logic to place orders based on demand forecast, current day inventory, planning calendar, order parameters and previously scheduled arrivals (ELT: expedited lead time, LT: lead time, PTF: planning time fence)}
  \label{fig:timeline}
\end{figure}

Material Requirements Planning (MRP) is at the heart of the process in today’s supply chain systems. MRP uses a deterministic logic that calculates the time and quantity of the orders that need to be placed, based on net material requirement, re-order parameters, current inventory, demand forecast and other inputs (Fig.~\ref{fig:timeline}). MRP computes material requirements based on demand forecast and schedules the arrivals as close as possible to the time they will be consumed. Re-order parameters such as Safety Stock Value (SSV) and Safety Time (ST) for the Safety Stock MRP system, or Reorder Point (RoP) for the Reorder Point MRP system are used as inputs to MRP to address the uncertainties in both supply and demand by providing buffers in the system.

MRP receives as input, 1) the bill of materials (BOM) which is a list of the materials, components and sub-assemblies required to make each product, 2) End Products  which is also called independent demand or level "0" on the BOM, 3) Quantity that is required at a given time, 4) Time at which the quantities are required to meet demand, 5) Inventory status records, which is a record of net materials available for use already in stock and materials on order from suppliers, and 6) Planning data which includes all the restraints and directions to produce such items as routing, labor and machine standards, quality and testing standards, lot sizing techniques (i.e. fixed lot size, lot-for-lot, economic order quantity), scrap percentages, and other inputs.

Demand forecast reflects the independent demand which is the demand for the \textit{end product} such as a computer or a bicycle. Dependent demand, on the other hand, is demand for \textit{component parts} or \textit{sub-assemblies}. For example, this would be the microchips in the computer, the wheels on the bicycle. Quantities for dependent demand are derived from independent demand using the Bill of Material (BOM). The planning system needs to consider not only the quantities of each of the component parts needed, but also the lead times needed to produce and receive the dependent demand items.

MRP can be expressed as an optimal control problem and the resulting optimal control problem can be solved by dynamic programming~\citep{zipkin2000foundations}. These dynamic programming methods (e.g., dynamic lot-sizing) were perceived to be too complex,  A few approximate heuristics (e.g., the Silver-Meal heuristic, 1973) were developed for the problem.  Such heuristics which are used in MRP systems to date, rely on gross simplifications. For example, Silver–Meal heuristic is a forward-looking method that requires determining the average cost per period as a function of the number of periods the current order is to span and stopping the computation when this function first increases. Well-known methods to find order quantities are Dynamic lot-sizing~\citep{wagner1958dynamic} introduced in 1958, Silver–Meal heuristic~\citep{silver1973heuristic} composed in 1973, and Least-unit-cost heuristic. As a result of such simplifying heuristics, the solution of the MRP in existing software are not optimal.

Commercially available MRP solutions suffer from three major drawbacks:

\begin{itemize}
  \item The dynamic programming solution to the optimal control formulation of MRP algorithms were perceived to be too complex. Hence, commercial MRP systems rely on grossly simplified approximate heuristics solutions to the problem which do not guarantee optimality of the MRP solution.

    \item MRP relies on input parameters to account for uncertainties associated with supply and demand functions that are set by end-users in an ad-hoc fashion, without rigorous analysis of uncertainty associated with historical supply and demand.

 \item Approximate solutions to MRP do not learn from nor adapt to the uncertainties and the time-varying and customer- or supplier-dependent operational constraints surrounding the manufacturing process.
\end{itemize}

Such limitations along with disparate, siloed data across multiple enterprise resource planning (ERP) systems resulting in stale and limited insights for the planners have led to overly conservative inventory estimates due to hard to quantify uncertainties.

Furthermore, a lack of end-to-end visibility and a system-level optimal decision making makes today's supply chains suffers form several issues. Delayed arrivals of purchase parts by supplier, causes production delay for downstream components and delay in production of components results in production delay for top-level products. Therefore, purchase part delays lead to overall trapped inventory and increased inventory holding cost. Most manufacturing processes suffer from poor forecasts. Demand forecasts end up being either under-forecasts which lead to stock-outs and poor OTIF, or over-forecasts which lead to excess inventory and hundreds of millions of dollars of capital locked up in unused inventory.


These challenges are presented as daunting problems in supply chain management and we propose to solve them with a stochastic inventory optimization (SIO) approach. In this article we will discuss a new formulation of inventory optimization problem as a constrained stochastic optimization problem in Section~\ref{section:methodology}. We discuss the data integration and model deployment platform that enabled solving this problem at large scale in Section~\ref{section:platform}. We present results showing inventory reductions across multiple systems at a global scale in Section~\ref{section:results}.

\nomenclature[A]{BOM}{Bill of materials}
\nomenclature[A]{ERP}{Enterprise resource planning}

\section{Related Work}
\label{section:related_work}

Solutions the MRP algorithm provides a heuristic approach for computing a sub-optimal but feasible solution for optimization problem, given SSV, ST, demand forecast, and other ordering parameters. To account for such limitations of MRP, new variations of MRP has been introduced. For example, \citet{ptak2011orlicky} introduced a new type of MRP called demand driven MRP (DDMRP). DDMRP leverages knowledge from theory of constraints, traditional MRP \& distribution resource planning (DRP), Six Sigma and lean manufacturing. It leverages MRP for planning, and kan-ban techniques for execution across multi-echelon supply chains, which means that it has the strengths of both but also the weaknesses of both and it remains a niche solution \citep{velasco2020applicability, smith2014demand}.

Since the multi-commodity distribution system design was first introduced by the work of \citet{inbook}, multiple optimization-based approaches have been proposed for the design of supply chain networks~\citep{AIKENS1985263, Geoffrion1995TwentyYO, VIDAL19971}. Most of these research, however, made the assumption that the supply chain are deterministic without considering the uncertainties from different sources such as demands and resource capacity. Supply chain disruption due to the ignorance of uncertain operating conditions will cause huge economic impact. A stochastic programming model was proposed by~\citet{SANTOSO200596} to solve the supply chain network design problem.

Since then, many attempts have been made to solve stochastic nature of supply chain network by leveraging mathematical methods. formulated A two-stage stochastic program targeting at an optimal operation plan was formulated by \citet{Ierapetritou1994NovelOA} and a decomposition-based optimization approach was proposed to solve the program. Further improvement on the approach was made by incorporating a MILP planning model~\citep{Gupta2000ATM} . A non-linear mixed integer formulation was proposed by \citet{Gupta2000ATM} with non-convex objective and constraint functions.

The main issue with the methods above is that they have oversimplified the supply chain in order to model it as a mathematical optimization, which has resulted in the omission of critical real-world operational complexities~\citep{agarwal2019multi}. Furthermore, it order to make the problem solvable additional relaxing assumptions have been made which cannot be made on real examples in real world operations. Simulation-based method was brought into the picture because of its capability of taking into account specific business restrictions, control variables and uncertainties of the problem~\citep{KOCHEL2005505}. However, a pure simulation-based approach won't lead to identification of best parameters to operate the system with, which could lead to optimal performance. As a result, a simulation-optimization approach has has been adopted recently to adequately model the complex system-level interactions and constraints and solve for optimal re-order parameters which leads to optimized inventory level and cost~\citep{10.5555/2700739.2700741}.

\citet{1172871} provided an overview of simulation optimization methods which were later adopted to solve the supply chain optimization problems. Given the involvement of multiple decision makers in the whole supply chain network. \citet{swaminathan1998modeling} proposed a multi-agent approach to represent critical supply chain components (e.g., retails, manufacturers) and simulate the supply chain based on their decisions. \citet{mele2006simulation} further extended agent-based  simulation-optimization approach by incorporating genetic algorithm. \citet{1166395} leveraged a multi-objective GA methodology with an existing supply chain simulator to optimize the system parameters. Besides directly optimize the inventory which may be limited by operation constraints, a simulation-optimization framework was proposed by \citet{JUNG20082570} to optimize the safety stock levels with consideration of production capacity.

The majority of methods mentioned above did not prioritize treatment of the stochastic nature of inventory above the optimization. For the ones that formulate the issue as stochastic programming problem, they lack the details in the modeling to properly represent the dynamics of inventory in supply chain. Those drawbacks make it difficult to productionize those solutions at a large scale.

\section{Methodology}
\label{section:methodology}

\subsection{Problem formulation}
\label{subsec:problem_formulation}

We cast the problem at hand as a stochastic constrained optimization problem and solve for the optimal set of re-order parameters that minimize inventory levels subject to a given service level constraint. The objective function in the optimization problem includes costs associated with holding inventory, ordering and transportation. This cost is optimized with respect to operational constraints. These operational constraints include planning calendar, planning time fence, (expedited) lead time, rounding values, minimum order quantity, customer-dependent operational constraints, complex material pricing (e.g., volume-based pricing). Another important constraint is often the service level constraint which is defined as expected probability of not hitting a stock-out during the next replenishment cycle.

The stochasticities, or uncertainties, in the constrained optimization problem come partly from the fact that there are often unexpected delays and shortages on the arriving orders. In addition, the miscellaneous and blocked material movements that cannot be anticipated ahead of time are also categorized as uncertainties. Finally, the optimization result is required to be frequently adjusted based on the latest estimation of the future demand, which is subject to change from one day to the next and is often different from the actual demand, leading to another source of uncertainty.

To address the constraints and the uncertainties mentioned above, we formulate the problem in a model predictive control (MPC) setting, which is a well-established framework in the process control discipline~\cite{garcia1989model, mayne2000constrained}. At each timestamp $t$ (also referred to as the sampling time), our problem can be formulated as the following optimization problem
\begin{equation}
\label{equ:mpc_formulation}
\begin{aligned}
    \min_{R^{:H}} \quad & \gamma_{p}\left[\sum_{i=0}^{H-1}\mathbf{x}^{t+i|t}\cdot \pi^{t+i}_{h} + \mathbf{a}^{t+i|t}\cdot\pi^{t+i}_{o,t}\right], \\
    \textrm{s.t.} \quad & \psi(\mathbf{X}^{:H}) - \text{SL}_{\text{min}} \geq 0, \\
    & \mathbf{x}^{t+i+1|t} = \mathbf{x}^{t+i|t} + \mathbf{a}^{t+i|t} - c_{\textrm{DF}}^{t, t+i} + \boldsymbol{\Theta} \\
    & \quad \forall i \in [0, \dots, H-1], \\
    & \mathbf{A}^{i+1:H} = \textrm{MRP}(\mathbf{x}^{t+i|t}, C_{\textrm{DF}}^{t+i}, r^{t+i|t}, \mathbf{A}^{i:H}) \\
    & \quad \forall i \in [0, \dots, H-1], \\
    & r^{t+i+1|t} = r^{t+i|t} \quad \forall i \in [0, \dots, H-1], \\
    & R^{:H} \geq \mathbf{0},
\end{aligned}
\end{equation}
where the variable $H$ denotes the length of the prediction horizon, which limits the length of the inventory simulation we perform at sampling time $t$. The superscript represents the timestamp, and the $t+i|t$ notation in the superscript depicts the value on timestamp $t+i$ predicted at the timestamp $t$. The variable $\mathbf{x}\in\mathbb{R}^{1}$ denotes the stochastic inventory level, and $\mathbf{a}\in\mathbb{R}^{1}$ depicts the stochastic expected order arrival quantity. Note that the expected order arrival quantity can be further split into the expected standard arrival (SA) and the expected expedited arrival (EA), which will be discussed in Section~\ref{subsec:mrp_sim}. The variable $\pi_{h}$ is the inventory holding cost, and $\pi_{o,t}$ represents both ordering and transportation costs. Note that objective function could take different forms depending on the business setting our customer has.

\nomenclature[V]{$\mathbf{x}$}{The stochastic inventory level for one day (scalars)}
\nomenclature[V]{$\mathbf{X}^{:\tau}$}{The 1d array of stochastic inventory level for a period of time, which starts at $t$ (inclusive) and ends at $t+\tau$ (exclusive)}
\nomenclature[V]{$\pi_{h}$}{Inventory holding cost}
\nomenclature[V]{$\pi_{o, t}$}{Inventory ordering and transportation cost}
\nomenclature[V]{$\gamma_{p}$}{Aggregation function that takes the $p$-th percentile of the uncertainty distribution}
\nomenclature[V]{$\psi$}{Function that calculates the service level}
\nomenclature[V]{$R^{:\tau}$}{The 1d array of reorder parameters for a period of time, which starts at $t$ (inclusive) and ends at $t+\tau$ (exclusive)}
\nomenclature[V]{$r$}{The reorder parameters for one day}
\nomenclature[V]{$\boldsymbol{\Theta}$}{Set of stochastic variables}
\nomenclature[V]{$c_{\textrm{DF}}^{t_o, t_f}$}{Demand forecast made at the date $t_o$ for the date $t_f$ (scalar)}
\nomenclature[V]{$C_{\textrm{DF}}^{\tau}$}{Demand forecast (1d array) made at the date $\tau$ for a length of predictive horizon $H$}
\nomenclature[V]{$C_{\textrm{DF, 2D}}$}{The 2-dimensional demand forecast matrix within the uncertainty sampling window}
\nomenclature[V]{$C_{\textrm{DF, 1D}}$}{Demand forecast aggregated to 1-dimension (from the off-diagonal elements in the 2d matrix) within the uncertainty sampling window}
\nomenclature[V]{$H$}{Predictive horizon}
\nomenclature[V]{$\mathbf{A}^{\tau:H}$}{The 1d array of the stochastic expected order arrivals starting at timestamp $t+\tau$ (inclusive) and ends at $t+H$ (exclusive)}
\nomenclature[V]{$\tilde{\mathbf{A}}^{:H}$}{The 1d array of the stochastic expected order arrival for the current predictive horizon predicted at a prior sampling time $\tilde{t}$. The array starts at $t$ (inclusive) and ends at $t+H$ (exclusive)}
\nomenclature[V]{$\mathbf{a}$}{The stochastic expected order arrival for one day (scalar)}

\nomenclature[A]{MRP}{Material requirements planning}
\nomenclature[A]{SSV}{Safety stock value}
\nomenclature[A]{ST}{Safety time}
\nomenclature[A]{SL}{Service level}
\nomenclature[A]{SIO}{Stochastic inventory optimization}

The notation $\gamma_{p}(\cdot)$ in the objective function of Eq.~\ref{equ:mpc_formulation} is the aggregation function that takes the specific percentile of a given uncertainty distribution, where the variable $p$ represents the percentile of the uncertainty distribution, (e.g., if $p=50\%$, $\gamma_{p}$ returns the median of the uncertainty distribution). We want to point out that aggregating with function $\gamma_{p}(\cdot)$ offers more flexibility than aggregating with the expectation (i.e., $\mathbb{E}[\cdot]$), although the latter is more often used as the as the aggregation function for stochastic optimization problems.

As mentioned in Section~\ref{section:introduction}, in all large scale operational settings where we intend to deploy our C3 SIO solution, the reordering of materials is done via legacy MRP systems. In order to respect that constraint and to ease the customer adaptation, we optimize for the reorder parameters, which are the inputs of the MRP system. Therefore, the decision variables in this formulation are denoted by $R^{:H}$, which represents the reorder parameters of the MRP system up to step $t+H$. To be specific, $R^{:H}$ is expressed as Eq.~\ref{equ:decision_variables} below
\begin{equation}
\label{equ:decision_variables}
    R^{:H} = [r^{t|t}, r^{t+1|t}, \dots, r^{t+H-1|t}] \in \mathbb{R}^{H\times 2}
\end{equation}
where the variable $r$ at a given timestamp is a single pair of reorder parameters for the MRP system. Note that our formulation would work for MRP systems with a variety of replenishment strategies. However, throughout the paper, we limit the scope of our discussion within the Safety Stock MRP system. Therefore, the reorder parameters here include the safety stock value (SSV) and the safety time (ST), and the variable $r^{t+i|t}$ can be represented by Eq.~\ref{equ:decision_variables_step} below
\begin{equation}
\label{equ:decision_variables_step}
    r^{t+i|t} = [\textrm{SSV}^{t+i|t}, \textrm{ST}^{t+i|t}] \in \mathbb{R}^{2}
\end{equation}

The variable $\mathbf{X}^{:H}$ in the first constraint of Eq.~\ref{equ:mpc_formulation} represents the vector of stochastic inventory level in the prediction horizon $H$, and is clearly defined in Eq.~\ref{equ:inventory_level} below:
\begin{equation}
\label{equ:inventory_level}
    \mathbf{X}^{:H} = [\mathbf{x}^{t|t}, \mathbf{x}^{t+1|t}, \dots, \mathbf{x}^{t+H-1|t}] \in \mathbb{R}^{H}
\end{equation}
where the bold lower-case notation $\mathbf{x}\in\mathbb{R}^{1}$ with superscripts represents the stochastic inventory level for a specific date, predicted at time $t$. The variable $\text{SL}_{\text{min}}$ represents the minimum target service level. The function $\psi(\cdot)$ calculates the service level within the prediction horizon. It can be expressed as Eq.~\ref{equ:service_level} below
\begin{equation}
\label{equ:service_level}
    \psi(\mathbf{X}^{:H}) = \gamma_{p_{\textrm{SL}}}\left[\frac{\sum_{i=1}^{H}\delta(\mathbf{x}^{t+i|t})}{H}\right]\in\mathbb{R}^{1},
\end{equation}
where the notation $\delta(\cdot)$ is a level set function, which is defined as
\begin{equation}
\label{equ:level_set}
    \delta(x) =
    \begin{cases}
        1, & x \geq 0\\
        0, & \textrm{otherwise}
    \end{cases}
\end{equation}
Note here the service level is computed as the percentage of days on which the inventory level is not negative. There are other definitions of the service level, such as fill rate (details in \cite{chu2015simulation}), which is also supported by this formulation. Also note that the aggregation functions $\gamma_{p}$ (in Eq.~\ref{equ:mpc_formulation}) and $\gamma_{p_{\textrm{SL}}}$ (in Eq.~\ref{equ:service_level}) are identical, but $p$ and $p_{\textrm{SL}}$ could take different values, depending on the risk profile of the customer. We use different notation here for clarity.

\nomenclature[V]{$p_{\textrm{SL}}$}{Service level percentile (SLP)}
\nomenclature[V]{$\delta$}{Level set function}

The second constraint in Eq.~\ref{equ:mpc_formulation} is a simple inventory updating logic, where $c_{\textrm{DF}}^{t, t+i}\in \mathbb{R}^{1}$ denotes the demand forecast for timestamp $t+i$ made at timestamp $t$, the variable $\boldsymbol{\Theta}$ denotes a collection of stochastic variables representing various sources of uncertainties. It includes the movement uncertainty $U_{MM}$, the supplier quantity uncertainty $U_{SQ}$, the supplier time uncertainty $U_{ST}$, and the demand forecast uncertainty $U_{DF}$, all of which will be discussed in detail in Section~\ref{subsec:uncertainty}.

The function notation $\text{MRP}(\cdot)$ in the third constraint in Eq.~\ref{equ:mpc_formulation} represents the constraints imposed by the legacy MRP system. It takes the inventory level up to the current timestamp $\mathbf{x}^{t+i|t}$, the demand forecasts at the current timestamp $C^{t+i}_{\textrm{DF}}$, the reorder parameters at the current timestamp $r^{t+i|t}$, and the expected order arrivals starting at the current timestamp $\mathbf{A}^{i:H}$, then it updates (overwrites) the expected order arrivals starting at the next timestamp to the end of the predictive horizon (i.e., $\mathbf{A}^{i+1:H}$).
%
The variable $C_{\textrm{DF}}^{t+i} \in \mathbb{R}^{H}$ is the demand forecast made at timestamp $t+i$, and it can be expressed as
\begin{equation}
\label{equ:demand_forecast_timestamp}
    C_{\textrm{DF}}^{t+i} = [c_{\textrm{DF}}^{t+i, t+i}, \dots, c_{\textrm{DF}}^{t+i, t+i+H-1}] \in \mathbb{R}^{H}.
\end{equation}
%
%
The notation $\mathbf{A}^{\tau:H}$ represents the stochastic expected order arrivals (include both standard arrivals and expedited arrivals) starting at timestamp $t+\tau$ ($\forall \tau \in [0, \dots, H-1]$) to the end of the predictive horizon (i.e., timestamp $t+H$). More specifically,
\begin{equation}
\label{equ:order_arrivals}
    \mathbf{A}^{\tau:H} = [\mathbf{a}^{t+\tau|t}, \dots, \mathbf{a}^{t+H-1|t}] \in \mathbb{R}^{H-\tau}.
\end{equation}
Note that the initial value for $\mathbf{A}^{i:H}$ in Eq.~\ref{equ:mpc_formulation} while $i=0$ (noted as $\tilde{\mathbf{A}}^{:H}$) comes from the the solution of Eq.~\ref{equ:mpc_formulation} on the previous sampling time $\tilde{t}$, which can be explicitly written as
\begin{equation}
\label{equ:order_arrivals_initial}
    \tilde{\mathbf{A}}^{:H} = [\mathbf{a}^{t|\tilde{t}}, \dots, \mathbf{a}^{t+H-1|\tilde{t}}] \in \mathbb{R}^{H}.
\end{equation}
More details on this MRP constraint can be found in Algorithm~\ref{algo:mrp_safety_stock}.

The fourth constraint in Eq.~\ref{equ:mpc_formulation} represents the limitation of the current MRP system, where only a single value can be specified for each of SSV and ST through out the entire prediction horizon. Note that the reorder parameters are not required to be identical when the sampling time $t$ evolves. For example, $r^{t+i|t}$ is not necessarily equal to $r^{t+i|t+j}$ for $j>0$ and $i>j$. The fifth constraint in Eq.~\ref{equ:mpc_formulation} indicates that the reorder parameters are always non-negative.

The formulation presented in Eq.~\ref{equ:mpc_formulation} will be solved iteratively for different sampling time $t$. The reorder parameters recommended at the sampling time $t$ will be preserved till it is updated at the next sampling time $t+j$. Note that the sampling time increment $j$ (also referred to as the optimization frequency) is not required to be 1. More details on the optimization frequency will be discussed in Section~\ref{subsec:detailed_model_treatements}. Also note that the reason we would expect to see non-zero results when we optimize the decision variables is the existence of the uncertainties (represented by $\boldsymbol{\Theta}$). Should the uncertainties in the formulation be completely removed, the decision variables will be all 0 since legacy MRP systems are capable of handling the deterministic inventory planning without any buffer specified.

\subsection{Algorithm overview}
Once the objective function and the constraints are configured (e.g., the hyper-parameters including $\gamma_{p}$ in Eq.~\ref{equ:service_level} determined), Eq.~\ref{equ:mpc_formulation} can be solved for at each sampling time. However, in practice, it is difficult for the customer to provide those hyper-parameters beforehand. Therefore, as part of the solution, we need to help the customer to configure the hyper-parameters that fit best for their business interest. In addition, MRP systems can be highly customized depending on the customer. There are pieces of logic such as order cancellation and order merging (i.e., flush window) that are difficult if possible at all to formulate as simple linear constraints for the optimization problem. This fact prevents us from directly feeding the optimization problem formulated in Eq.~\ref{equ:mpc_formulation} in a traditional mixed-integer linear programming (MILP) solver.

To handle the challenges mentioned above, the C3 AI Inventory Optimization application introduces an AI Stochastic Inventory Optimization Algorithm (``the algorithm"). The algorithm learns a material plant's historical uncertainties in material movements - supply, demand, and other - over a historical time period (``training period") to select the best hyper-parameters that are feasible to the customer. Then, the algorithm perturbs these learned material movement uncertainties into the material planning simulations in the future (i.e., ``validation period", and ``live production"). Based on those simulations and uncertainty realizations in the future, and with the help of the selected hyper-parameters, the algorithm recommends optimized reorder parameters (i.e., SSV and ST for the Safety Stock MRP system) that minimize inventory costs while meeting the target service level constraint.

\begin{figure*}[!htp]
  \centering
  \includegraphics[width=\linewidth]{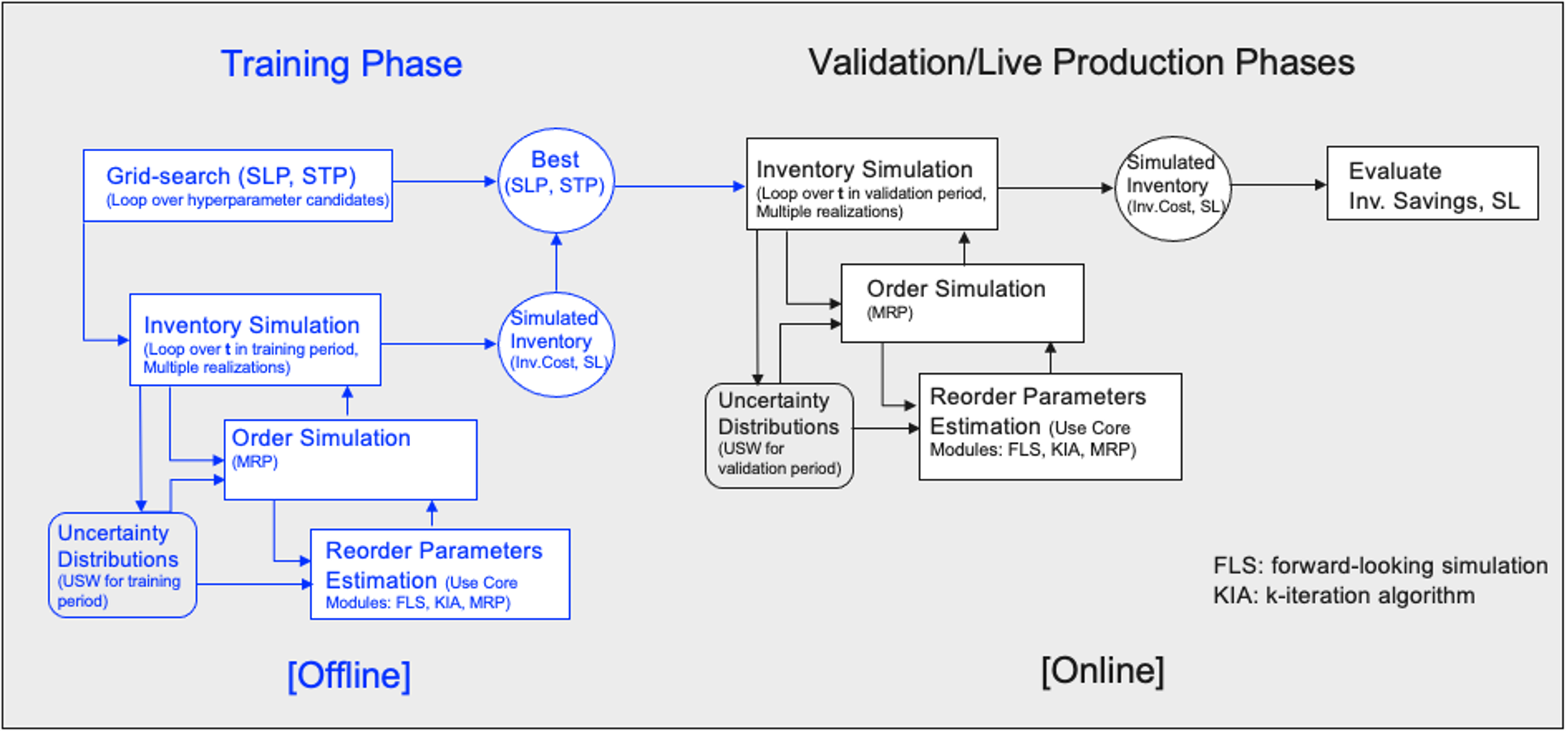}
  \caption{Workflow of the C3 AI Stochastic Inventory Optimization algorithm}
  \label{fig:algorithm_workflow}
\end{figure*}

Fig.~\ref{fig:algorithm_workflow} demonstrates the general workflow of the algorithm. The algorithm has two phases: a training (offline) phase and a validation or live production (online) phase. The training phase happens on a period of time further in the history while the validation phase is on a more recent time period. The term ``simulation period" in this paper refers to the timestamps (i.e., days) either in the training or in the validation phase depending on the context. Note that the live production is a special case of the validation phase where the optimization is performed at ``today" (the last date we have the data available). Therefore, the reorder parameters are only optimized for one timestamp in the live production phase.

\nomenclature[A]{FLS}{Forward-looking simulation}
\nomenclature[A]{KIA}{The k-iteration algorithm}

We start the discussion of Fig.~\ref{fig:algorithm_workflow} from the validation phase, where the optimization problem formulated in Eq.~\ref{equ:mpc_formulation} is solved repeatedly at each sampling time $t$ to generate the final simulated inventory for the validation period. The inventory simulation of the validation phase has three key components (steps): (1) uncertainty distribution step, (2) reorder parameters estimation step, and (3) order simulation step. During the uncertainty distribution step, all the uncertainties of the model inputs are collected within the uncertainty sampling window (USW). Those uncertainties are saved for the sampling procedure that happens in the reorder parameters estimation step. In the reorder parameters estimation step, the MRP simulations are performed repeatedly starting at each sampling time $t$ to infer the best reorder parameters. The inference of the re-order parameters is accomplished by the ``k-iteration" process (with the hyper-parameters determined by the training phase), which will be discussed in Section~\ref{subsec:k_iteration}. The reorder parameters will be retained till they are overwritten by the inference made at the next sampling time. The algorithm is designed with this overwriting process to effectively incorporate the new information that becomes available at each sampling time $t$.

Subsequently, in the order simulation step, the future order placements will be generated based on a series of reorder parameters. Finally, the algorithm incorporates the uncertainty distributions again to generate multiple realizations of the simulated inventory. The inventory cost of each realizations can then be computed accordingly. Inventory savings can be computed by comparing the actual cost of the inventory with the median (or any other aggregation approach configurable by the user) of the realizations during the validation period.

The live production (or, operation) is essentially identical to the validation phase. The only distinction is that there is only one timestamp (day) for reorder parameter recommendation in the live production. Also, since we will not have the actual data (e.g., actual inventory, actual arrivals, actual consumption, etc.) for the live production, there is no need to compare the simulated results with the actual inventory.

The purpose of the training phase is to determine the hyper-parameters, which helps the inferring of the reorder parameters at each sampling time $t$ in the validation phase. The training phase (shown in blue in Fig.~\ref{fig:algorithm_workflow}) is essentially a grid search process repeating the validation procedure multiple times in the training period, each time with a different combination of hyper-parameters. Then the best hyper-parameter is determined based on the inventory savings and service level achieved by each combination. Note that the grid search approach is chosen among all optimization options to ensure the stability of the algorithm, which will be discussed in Section~\ref{subsec:opt_hyperparams}.

To summarize, the deployment of the algorithm includes three key steps:
\begin{enumerate}
    \item Train the model: the model learns from material movements and associated uncertainties over the historical training period and optimizes the hyper-parameters.
    \item Validate the model: the trained model (with optimized hyper-parameters) is validated in a production-like operation over a period of recent history to evaluate real performance (i.e., inventory cost and service level).
    \item  Operate the model (live production): The model is run in production, simulating uncertainties into the future to provide recommendations on the optimal reorder parameters.
    The model's performance is monitored; it may be automatically or manually retrained as needed.
\end{enumerate}

In the following parts of this section, we will cover topics including the MRP simulation, the modeling of uncertainties, optimization of the hyper-parameters, the forward-looking simulation, the ``k-iteration" algorithm, and other algorithm treatments.

\subsubsection{Inventory simulation with MRP}\label{subsec:mrp_sim}
As introduced previously, the MRP system sits at the core of C3 AI SIO Algorithm. For any given sampling time $t$, the MRP takes the reorder parameters along with other inputs including demand forecast ($C_{\textrm{DF}}^{t+i}$), lead time (LT), expedited lead time (ELT), planning calendar (PC), item volume pricing (IVP), etc. within a horizon (corresponds to the predictive horizon $H$ in Eq.~\ref{equ:mpc_formulation}, then outputs the order placement in that horizon.

As is explained in Section~\ref{subsec:problem_formulation}, the optimization framework introduced by C3 AI SIO could support MRP systems with a variety of replenishment strategies. However, we limit the scope of our discussion in this paper to the Safety Stock MRP where the reorder parameters are safety stock value (SSV) and safety time (ST).

\nomenclature[A]{LT}{Lead time}
\nomenclature[A]{ELT}{Expedited lead time}
\nomenclature[A]{RV}{Rounding value}
\nomenclature[A]{MO}{Minimum order quantity}
\nomenclature[A]{MO}{Minimum order quantity}
\nomenclature[A]{PTF}{Planning time fence}
\nomenclature[A]{SA}{Expected standard arrival}
\nomenclature[A]{EA}{Expected expedited arrival}

\paragraph{Safety Stock MRP}\mbox{}\\

The Safety Stock MRP system places orders for a window in the future based on the current inventory position, existing orders and demand forecast. The MRP projects the inventory position for every day in this future window. Orders are placed to arrive when the projected inventory drops below a predefined safety stock value, subject to reorder constraints such as lead time, factory calendar, planning calendar, minimum order quantity and rounding values, etc.

Algorithm~\ref{algo:mrp_safety_stock} provides a simplified view of the Safety Stock MRP system, where the key pieces of logic are described in details. However, some auxiliary pieces logic including item volume pricing (IVP), lot size in days (LS), planning calendar (PC), order flushing (or, merging) are intentionally left out in this algorithm to improve the readability of the algorithm. Order cancellation logic, on the other hand, is explained in the next sub-section (Algorithm~\ref{algo:order_cancellation}).

Input data to Algorithm~\ref{algo:mrp_safety_stock} include starting inventory $\bx^t$, demand forecast $C_{\textrm{DF}}^{t}$, the expected arrivals of existing standard orders $\bA^{t+i}$ and existing expedited orders $\bEA^{t+i}$, where $i=0,...,H$, and $H$ is the horizon. The algorithm simulates inventory for the horizon $\bxh^{t+i}, \hspace{2px} i=0,...,H$. Notice that hat accent is used to denote inventory position within MRP. When $ELT<=i<LT$, where $ELT$ and $LT$ are expedited lead time and lead time, respectively, the algorithm places new expedited orders if projected inventory $\bxt^{t}$ is about to go below $0$. When $LT<=i<H$, the algorithm places new standard orders if projected inventory is about to go below safety stock value $SSV$. The quantity of new standard orders are modified to comply with minimum order quantity $MO$ and rounding value $RV$.

\begin{algorithm}
\caption{Simplified Safety Stock MRP}
\label{algo:mrp_safety_stock}
\hspace*{\algorithmicindent} \textbf{Input:} starting inventory $\bx^t$, demand forecast $C_{\textrm{DF}}^{t}$, expected standard order arrivals $\bA^{t+i}$, expected expedited order arrivals $\bEA^{t+i}$, where $i=0...H$, and $H$ is the planning time horizon, lead time $LT$, expedited lead time $ELT$, safety stock value $SSV$, minimum order quantity $MO$, rounding value $RV$\\
\hspace*{\algorithmicindent} \textbf{Output:} updated expected arrivals of standard orders $\bA^{t+i}$ and expedited orders $\bEA^{t+i}$, $i=0...H$
\\
\begin{algorithmic}[1]
\State Apply order cancellation on expected standard order arrivals ($\bA$) based on Algorithm~\ref{algo:order_cancellation}
\State $\bxh^t = \bx^t$
\\
\For{$i = 0, \dots, ELT - 1$}
    \State $\bxh^{t+i+1} = \bxh^{t+i} + \bA^{t+i} + \bEA^{t+i} - c_{\textrm{DF}}^{t, t+i}$
\EndFor
\\
\For{$i = ELT, \dots, LT - 1$}
    \State $\bxt^{t+i+1} = \bxh^{t+i} + \bA^{t+i} + \bEA^{t+i} - c_{\textrm{DF}}^{t, t+i}$
    \If{$\bxt^{t+i+1} < 0$}
        \State $\bEA^{t+i} \mathrel{{+}{=}} -\bxt^{t+i+1}$
    \EndIf
    \State $\bxh^{t+i+1} = \bxh^{t+i} + \bA^{t+i} + \bEA^{t+i} - c_{\textrm{DF}}^{t, t+i}$
\EndFor
\\
\For{$i = LT, \dots, H-1$}
    \State $\bxt^{t+i+1} = \bxh^{t+i} + \bA^{t+i} + \bEA^{t+i} - c_{\textrm{DF}}^{t, t+i}$
    \If{$\bxt^{t+i+1} < SSV$}
        \State $\Delta x = SSV - \bxt^{t+i+1}$
        \State $\bA^{t+i} \mathrel{{+}{=}} \max(\lceil\frac{\Delta x - MO}{RV}\rceil, 0) \times RV + MO$
    \EndIf

    \State $\bxh^{t+i+1} = \bxh^{t+i} + \bA^{t+i} + \bEA^{t+i} - c_{\textrm{DF}}^{t, t+i}$
\EndFor
\\
\State \Return{$\bA^{t+i}$, $\bEA^{t+i} \hspace{4px} i=0...H$}
\end{algorithmic}
\end{algorithm}

\paragraph{Order cancellation}\label{subsec:oc_sim}\mbox{}\\
Orders sent to vendors can be cancelled based on a cancellation window (e.g., 30 days prior to arrival). Usually the orders can be cancelled when the orders’ expected arrival date is later than the planning-time-fence ($PTF$) day, counting from the current day. Within the MRP module, we are enabled to cancel orders when the $PTF$ is greater than the lead time ($LT$). However, in the case that $PTF$ being smaller than $LT$, MRP does not cancel the orders that will be arriving in between $[PTF, LT]$ if more quantity is ordered than needed.

We developed a pre-processing algorithm (Algorithm~\ref{algo:order_cancellation}) for MRP to handle the cancellation of the excessive orders. The algorithm modifies the orders that were placed previously but are expected to arrive in between $[PTF, LT]$. Rather than iterating incrementally from $PTF$ to $LT$, it iterates in the backward order from $LT$ to $PTF$ to make sure the cancelled order does not impact future inventory levels. Within each iterated date, it computes the surplus of inventory, which is denoted by $\tilde{x}^{t+k}$ in the algorithm. If the surplus is greater than zero, the expected arrivals for that date would be reduced to offset the surplus.

\begin{algorithm}
\caption{The order cancellation logic}
\label{algo:order_cancellation}
\hspace*{\algorithmicindent} \textbf{Input:} inventory $x^{t}$ of the day $t$, expected standard order arrivals $\bA^{t+i}$, expedited order arrivals $\bEA^{t+i}$ and demand forecast $c_{DF}^{t,t+i}$ where $i = 0...H$, and $H$ is the planning time horizon for a single round of MRP simulation, lead time $LT$, planning time fence $PTF$, and safety stock value $SSV$\\
\hspace*{\algorithmicindent} \textbf{Output:} modified expected order arrivals $\bA^{t+i}$ where $i = 0...H$
\begin{algorithmic}[1]
\State Set $k = LT$
\While{$k > PTF$}
    \State $\tilde{x}^{t+k} = x^{t} + \sum_{i=0}^{i=k}\bA^{t+i} + \sum_{i=0}^{i=k}\bEA^{t+i} - \sum_{i=0}^{i=k}c_{DF}^{t,t+i} - SSV$
    \If{$\tilde{x}^{t+k} \leq 0$}
        \State $k = k-1$
        \State continue to next iteration
    \EndIf

    \If{$\tilde{x}^{t+k} > \bA^{t+k}$}
        \State Set $\bA^{t+k} := 0$
    \Else
        \State Set $\bA^{t+k} := \bA^{t+k} - \tilde{x}^{t+k}$
    \EndIf
    \State $k = k-1$
\EndWhile
\State \Return{$\bA$}
\end{algorithmic}
\end{algorithm}


\subsubsection{Modeling of the uncertainties}
\label{subsec:uncertainty}
In our algorithm, we consider four major sources of uncertainties, including demand forecast uncertainty, material movement uncertainty, supplier quantity uncertainty and supplier time uncertainty.

Depending on the accuracy of the demand forecast model, the predicted demand can be quite different from the actual demand/consumption. Demand forecast uncertainty models the distribution of the difference between demand forecast and actual consumption.

There are two types of material movements, miscellaneous movements, and blocked movements. These are typically not considered in traditional MRPs and are not part of the amount represented in the demand forecast. Miscellaneous movements happen due to movement between facilities or customer returns. Blocked movements happen when parts arrived but cannot be consumed (e.g., due to a quality issue). As these material movements are generally unpredictable, they are treated as sources of uncertainty. Material movement uncertainty models the distribution of the sum of miscellaneous movements and blocked movements for a given SKU.

Supplier quantity uncertainty models the distribution of the difference between the planned quantity in the purchase orders and the actual quantity received for a given SKU. Supplier time uncertainty models the distribution of the time difference of the planned delivery date and the actual delivery date of orders for a given SKU.

In reality, different suppliers can have different behaviors. One supplier may be more likely to delay and reduce quantity than another supplier. However, we usually don't have enough data points to accurately model the behavior of each supplier. Therefore, we average out different behaviors of all suppliers for a given SKU into one distribution.

We estimate the distribution of each source of uncertainty using an empirical distribution. Each empirical distribution is represented by a list of historical data points. We assume that all historical data points are equally likely, then we can sample the empirical distribution from the list of historical data points. The list of historical data points is retrieved from relevant quantities within a certain uncertainty sampling window (USW).

The USW for training is the entire training period. The USW during validation is a growing window with a fixed start date, where the initial length of the USW $L_\textrm{USW}^0$ is given by
\begin{equation}
\label{equ:usw-length}
    L^\textrm{USW}_0 =\max(L_\textrm{min}^\text{USW}, LT+b^{\textrm{USW}}),
\end{equation}
where $LT$ represents lead time, $b^{\textrm{USW}}$ represents the buffer (or, the offset) of the USW length and $L_\textrm{min}^\text{USW}$ denotes the minimum initial length of USW.

The empirical distribution for material movement uncertainty is represented by $U_\text{MM}$, which is a list of data points in the uncertainty sampling window, where each data point is the sum of miscellaneous movement and blocked movement at that time step. To be more specific,

\begin{equation}
\label{equ:mm-uncertainty}
    U_\text{MM} =[..., c_\text{MIM}^{\tau} + c_\text{BM} ^{\tau}, ...] \in \mathbb{R}^{L^{\textrm{USW}}},  \forall \tau \in \text{USW},
\end{equation}
where $L^\text{USW}$ is the length of the uncertainty sampling window, $c_\text{MIM}^{\tau}\in\mathbb{R}^{1}$ and $c_\text{BM} ^{\tau}\in\mathbb{R}^{1}$ denotes the value of consumption due to miscellaneous movement and blocked movement at time $\tau$, respectively. The variables $c_\text{MIM}^{\tau}$ and $c_\text{BM} ^{\tau}$ are positive when materials move in and negative when materials move out.

We use $U_\text{SQ}$ to denote the list of historical data points that represents the empirical distribution of the supplier quantity uncertainty. Each data point in $U_\text{SQ}$ is the difference between planned quantity and actual received quantity of an order that is planned to be delivered in the uncertainty sampling window. It can be expressed as

\begin{equation}
\label{equ:squ}
\begin{split}
    U_\text{SQ} = & [..., \min(s_\text{A}^o - s_\text{P
}^o, 0),...] \in \mathbb{R}^{L^{\textrm{USW}}}, \\
& \forall o \in \text{Orders planned to arrive in USW} ,
\end{split}
\end{equation}
where $s_\text{A}^o\in\mathbb{R}^{1}$ and $s_\text{P}^o\in\mathbb{R}^{1}$ denote the actual and planned supplied quantity from order $o$, respectively. To make sure our algorithm is conservative in nature, when the actual supplied quantity $S_\text{A}^o$ is larger than the planned supplied quantity $S_\text{P}^o$, we do not take this favorable condition into consideration when quantifying the uncertainty. Therefore, we will replace the positive number with 0 in the empirical distribution.

We use $U_\text{ST}$ to denote the list of historical data points that represents the empirical distribution of the supplier time uncertainty. Each data point in $U_\text{ST}$ is the difference (in number of days) between planned delivery date and actual delivery date of an order that is planned to be delivered in the uncertainty sampling window.

\begin{equation}
\label{equ:stu}
\begin{split}
    U_\text{ST} = & [..., \max(t_\text{A}^o - t_\text{P}^o, 0),...] \in \mathbb{R}^{L^{\textrm{USW}}}, \\
    & \forall o \in \text{Orders planned to delivered in USW},
\end{split}
\end{equation}
where $t_\text{A}^o\in\mathbb{R}^{1}$ and $t_\text{P}^o\in\mathbb{R}^{1}$ denote the actual and planned delivery time (in days) from order $o$, respectively. When the actual delivery date $t_\text{A}^o$ is earlier than the planned delivery date $t_\text{P}^o$, again, this is a favorable condition that we do not take into consideration. Therefore, we will instead put a value of 0 into the empirical distribution.

To construct the demand forecast uncertainty $U_\text{DF}$, we need to compress the 2-dimensional demand forecast expressed in Eq.~\ref{equ:demand_forecast}
\begin{equation}
\label{equ:demand_forecast}
    C_\text{DF, 2D} = [..., C_{\textrm{DF}}^{\tau}, ...] \in \mathbb{R}^{L^{\textrm{USW}}\times H}, \forall \tau \in \text{USW}.
\end{equation}
into a 1-dimensional array before subtracting it by the actual consumption. To be specific, following Eq.~\ref{equ:demand_forecast_timestamp}, let $c^{t_o, t_f}_\text{DF}\in\mathbb{R}^{1}$ be the demand forecast (planned consumption) for time $t_f$ made at time $t_o$. For each time $\tau$ in the uncertainty sampling window, extract the demand forecast value made for time $\tau$ on time $\tau - \textrm{LT} - \textrm{ST}$ to get a 1-dimensional representation of the demand forecast, where $\textrm{LT}$ and $\textrm{ST}$ stand for lead time and safety time, respectively. The 1-dimensional representation can be expressed as

\begin{equation}
\label{equ:demand_forecast_1d}
    C_\text{DF, 1D} = [..., c^{\tau - \textrm{LT} -\textrm{ST}, \tau}_\text{DF}, ...] \in \mathbb{R}^{L^{\textrm{USW}}}, \forall \tau \in \text{USW}.
\end{equation}
Here we use demand forecast made on lead time plus safety time days before as we believe this most accurately reflects the actual demand forecast uncertainty. Note that $C_\text{DF, 1D}$ effectively comes from the off-diagonal terms of the matrix $C_\text{DF, 2D}$.

In practice, demand forecast oftentimes comes in on a weekly basis or monthly basis. To convert the demand forecast into a daily one, so that it is comparable to the actual consumption, we apply a convolution step on $C_\text{DF, 1D}$ that acts as a moving average smoother. Specifically, when demand forecast comes in once everyday $d$, we convolve $C_\text{DF, 1D}$ with a $d$-dimensional constant vector $\mathbf{a} = [\dfrac{1}{d},...,\dfrac{1}{d}]$ to obtain the smoothed 1-dimensional demand forecast $C_\text{DF, 1D}^*$.
\begin{equation}
\label{equ:demand_forecast_convolve}
    C_\text{DF, 1D}^{i,*} = \dfrac{1}{d}\sum_{j=i-d/2}^{i+d/2}C_\text{DF, 1D}^j,
\end{equation}
where $C_\text{DF, 1D}^j$ represents the $j$-th element in $C_\text{DF, 1D}$ defined in Eq.~\ref{equ:demand_forecast_1d}, and $C_\text{DF, 1D}^{i,*}$ denotes the $i$-th element in the smoothed 1-dimensional demand forecast $C_\text{DF, 1D}^*$.

Let
\begin{equation}
    C_\text{A} = [..., c^{\tau}_\text{A}, ...] \in \mathbb{R}^{L^{\textrm{USW}}}, \forall \tau \in \text{USW}
\end{equation}
be the list of actual consumption in the uncertainty sampling window. We compute the demand forecast uncertainty with
\begin{equation}\label{equ:dfu}
    U_\text{DF} = C_\text{DF, 1D}^* - C_\text{A} \in \mathbb{R}^{L^{\textrm{USW}}}.
\end{equation}

Note that in practice, we would encounter outliers in $C_{\textrm{DF}}$, so pre-processing steps (discussed in Section~\ref{subsec:detailed_model_treatements}) are necessary before computing $U_{\textrm{DF}}$ or consuming $C_{\textrm{DF}}$ in the algorithm.

\nomenclature[V]{$c_{\textrm{A}}^{\tau}$}{Actual consumption at timestamp $\tau$}
\nomenclature[V]{$U_{\textrm{ST}}$}{The 1d array to represent the supplier time uncertainty within the uncertainty sampling window}
\nomenclature[V]{$U_{\textrm{SQ}}$}{The 1d array to represent the supplier quantity uncertainty within the uncertainty sampling window}
\nomenclature[V]{$U_{\textrm{MM}}$}{The 1d array to represent the material movement uncertainty within the uncertainty sampling window}
\nomenclature[V]{$U_{\textrm{DF}}$}{The 1d array to represent the demand forecast uncertainty within the uncertainty sampling window}
\nomenclature[A]{USW}{Uncertainty sampling window}
\nomenclature[V]{$N_r$}{Number of uncertainty realizations in the forward-looking simulation procedure}
\nomenclature[V]{$b^{\textrm{USW}}$}{The buffer for the length of the uncertainty sampling window}
\nomenclature[V]{$L_{\textrm{min}}^{\textrm{USW}}$}{The minimum length of the uncertainty sampling window}

\subsubsection{Optimization of the hyper-parameters}
\label{subsec:opt_hyperparams}

\nomenclature[A]{SLP}{Service level percentile}
\nomenclature[A]{STP}{Safety time percentile}

As mentioned earlier in the algorithm overview, we will help the customers to decide the risk profile that fits best to their business interest. The risk profile can be quantified by hyper-parameters including $p_{\textrm{SL}}$ for $\gamma_{p_{\textrm{SL}}}$ in Eq.~\ref{equ:service_level}. This hyper-parameter is also referred to as service level percentile, or SLP. As will be discussed later, in the context of Safety Stock MRP system, the SLP directly affects the choice of SSV in the C3 SIO algorithm, but it does not have any impact on the choice of ST. To help with quantifying customer's risk profile on the time delay and recommending ST, we introduce another hyper-parameter, safety time percentile (or, STP). Although those two variables are called ``hyper-parameters", they are analogous to the ``parameters" in a regular machine learning model, which is optimized during the training process and directly used in the model for the test period. They are referred to as ``hyper-parameters" in C3 SIO algorithm to distinguish from other model parameters that are not optimized during the training process.

To distinguish SLP from the service level (SL), we note that SL describes the performance of a single simulated inventory curve. However, SLP describes the performance of multiple simulated inventory curves, and it is not a replacement but a complement of SL. SLP describes the percentage of inventory realizations that meets the minimum target service level (i.e., $\textrm{SL}_{\textrm{min}}$ in Eq.~\ref{equ:service_level}). For example, if the minimum target service level is 96\% and the total number of uncertainty realization of simulated inventory is 100, which means there will be 100 simulated inventory curves. A service level percentile of 50\% means that 50 out of the total 100 realizations will have a service level above 96\%.

\begin{figure}[!htp]
  \centering
  \includegraphics[width=\linewidth]{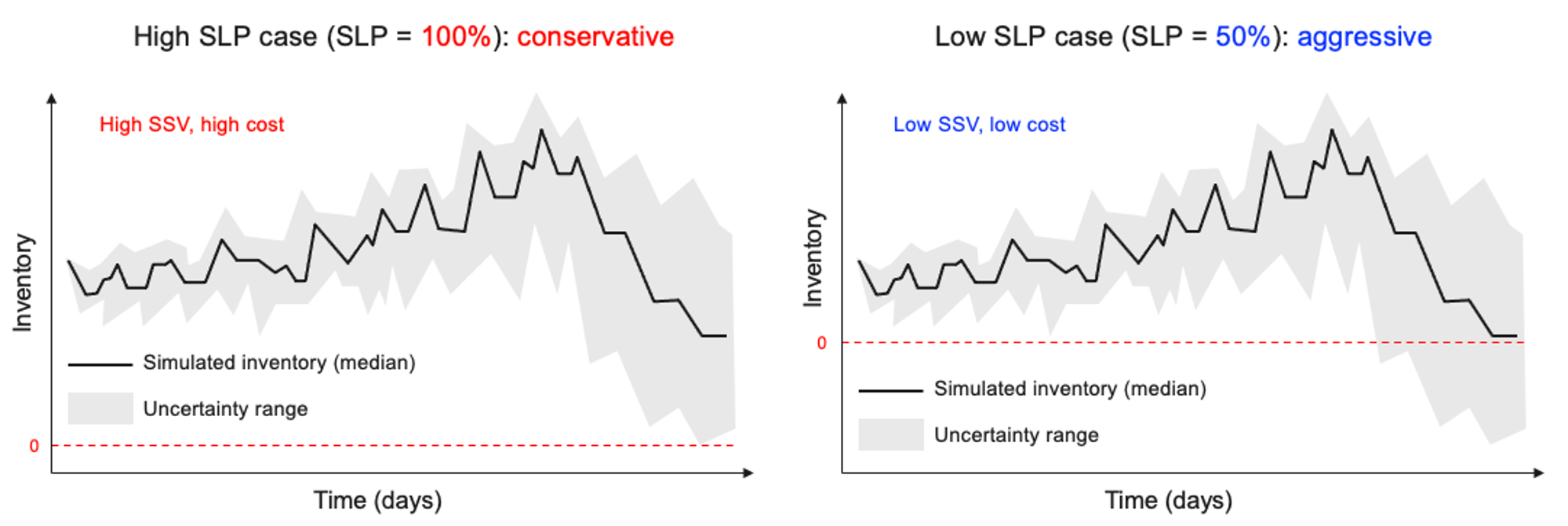}
  \caption{Comparison of high versus low service level percentiles}
  \label{fig:slp_illustration}
\end{figure}

SLP will help the algorithm to recommend SSV at each sampling time during the ``k-iteration" process, which will be described in detail in Section~\ref{subsec:k_iteration}. SLP decides how conservative or aggressive the SSV choice will be. Fig.~\ref{fig:slp_illustration} provides a visual illustration on the discrepancies between a high SLP and a low SLP with the grey areas represent the uncertainty range of the simulated inventory. Note that we are using a minimum target service level of 100\% while discussing about SLP in Fig.~\ref{fig:slp_illustration}. A high SLP (e.g., 100\% as is shown on the left of the figure) puts the entire uncertainty range of the simulated inventory above 0 (i.e., for 100 out of 100 uncertainty realizations, each of the realization has a service level of 100\%), while a low SLP (e.g., 50\% on the right side of the figure) only puts the median curve above 0. The low SLP is a more aggressive strategy, which means the inventory will have a higher chance of going below the minimum target service level (high risk for stock-out), but the inventory cost will be much lower (high return for inventory savings). For now, we pick the SLP in the training based on the inventory cost of the median realizations (curves in black, which means we set $p=50\%$ for $\gamma_{p}$ in Eq.~\ref{equ:mpc_formulation}), but it is configurable by the customer.

STP, on the other hand, describes the percentile in the supplier time uncertainty distribution. The strategy used in the Safety Stock MRP system to de-risk the supplier delay is to place the purchase order earlier. The safety time is the variable that determines how early the order should be placed compared to the default order placement date in the MRP output. Therefore, the choice of the most suitable safety time depends on the supplier time uncertainty profile, which could vary over the simulation period. A fixed choice of STP could map to different safety time while the supplier time uncertainty profile changes over time.

\begin{figure}[!htp]
  \centering
  \includegraphics[width=\linewidth]{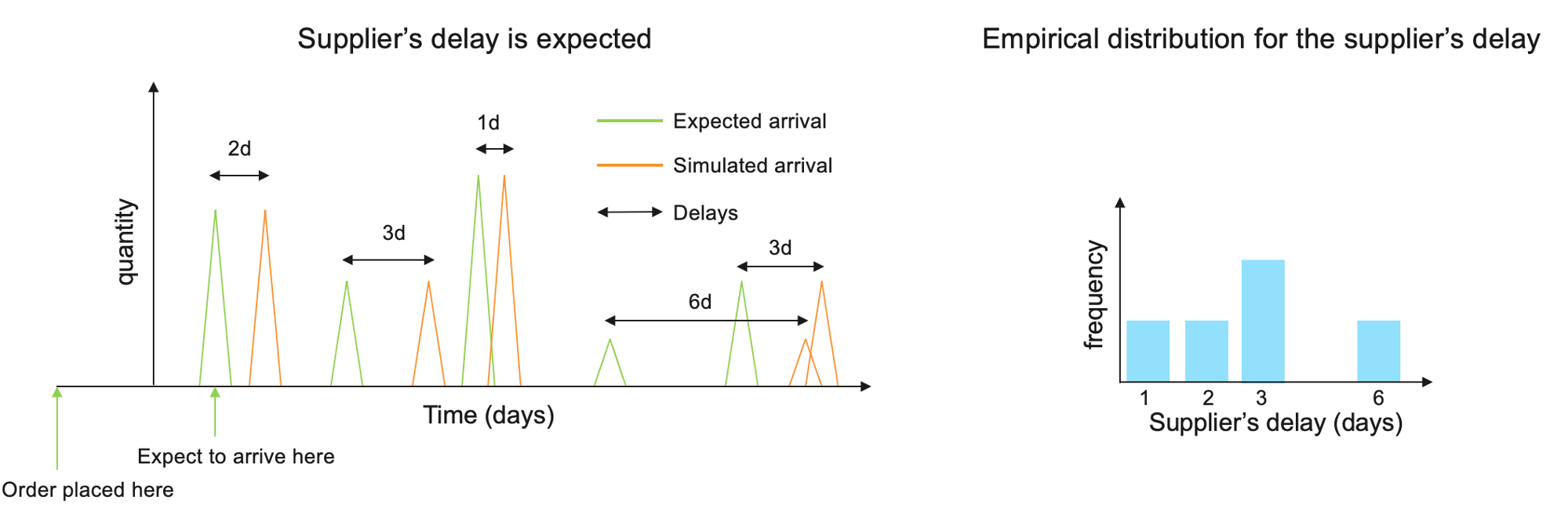}
  \caption{Safety time percentile is the percentile on the supplier time uncertainty}
  \label{fig:stp_illustration}
\end{figure}

Fig.~\ref{fig:stp_illustration} demonstrates what distribution are we looking at while deciding the safety time percentile. The plot on the left shows the delay of each arriving order in the USW (i.e., the delays are represented by the differences between the green peaks and the orange peaks). Those delays are translated into an empirical distribution of supplier time uncertainty as is shown in the right plot. Safety time will then be determined based on the STP and the empirical distribution. For example, a STP of 50\% would translate into a safety time of 3 days with the empirical distribution provided in Fig.~\ref{fig:stp_illustration}, since the 50th percentile of the empirical distribution is 3 days, and a STP of 100\% would translate into safety time of 6 days. Similar to the choice of SLP, a high STP represents a conservative strategy since it would make MRP to place all the orders so early that even the longest delay would be prevented. This STP would result in a high service level, however, the inventory level will be piling up due to the early placement of orders and the inventory cost will be high. On the other side, a low STP is a aggressive strategy since it will result in a smaller safety time. Therefore, the inventory is prone to stock-out, but the inventory cost will be lower.

\begin{figure}[!htp]
  \centering
  \includegraphics[width=\linewidth]{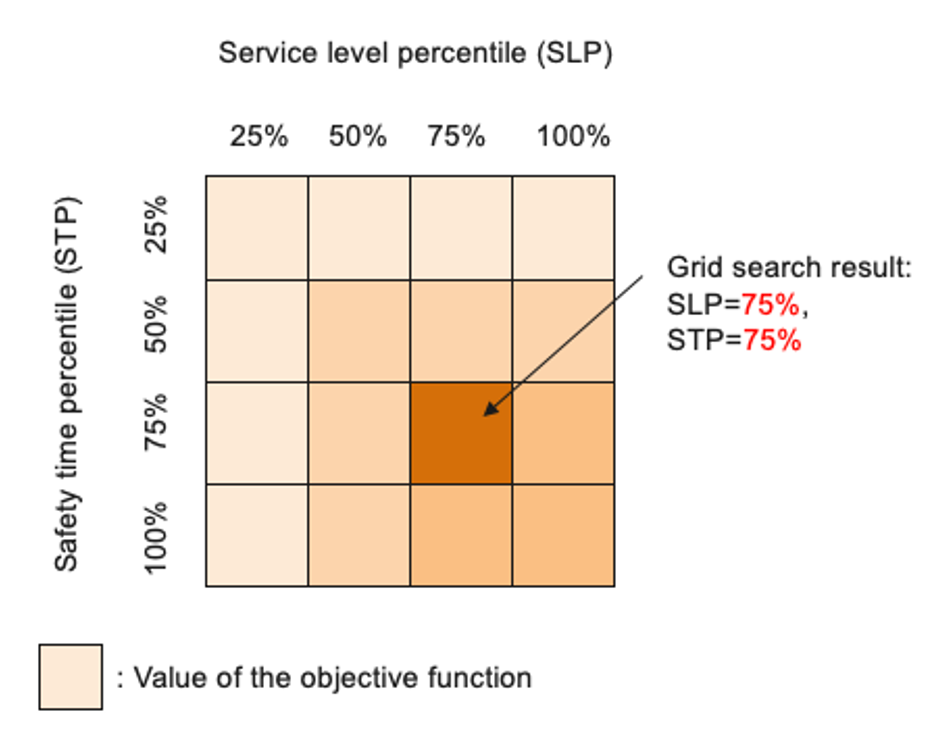}
  \caption{Grid search of hyper-parameters}
  \label{fig:grid_search}
\end{figure}

In the training period color-coded in blue in Fig.~\ref{fig:algorithm_workflow}, the best combination of SLP and STP is selected by the grid search approach (illustrated in Fig.~\ref{fig:grid_search}) so that the median of the simulated inventory realizations would have the lowest inventory cost while its service level is above the target service level. Note that we apply the objective function and the constraint on the median of the realizations as a design choice, which is configurable to fit the user's preferences.

\begin{figure}[!htp]
  \centering
  \includegraphics[width=\linewidth]{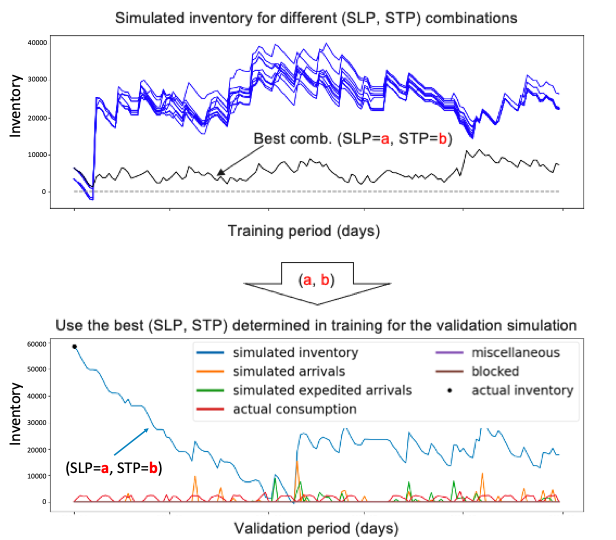}
  \caption{Selection and consumption of the hyper-parameters}
  \label{fig:grid_search_train_test}
\end{figure}

Fig.~\ref{fig:grid_search_train_test} provides a holistic picture of how the hyper-parameters are selected and consumed in the algorithm. During the training period, the median realization of the simulations for all the combination of SLP and STP are generated, as is shown in the top part of the figure. Among all the combinations, one combination with SLP equal to $a$ and STP equal to $b$ gives the lowest inventory cost while satisfies the constraint on the service level. Therefore, this combination of $(a,b)$ is selected as the best hyper-parameters. This hyper-parameter pair is then used in the validation period. In the validation period, the SLP and STP are translated to SSV and ST recommendations at each timestamp of the simulation. Then the inventory is simulated based on those recommendations to generate simulated inventory curve shown as the light blue curve in the bottom figure.

\subsubsection{Forward-looking simulations}
\label{subsec:forward_looking}

After selecting the hyper-parameters, we enter the reorder parameter estimation step. To find the optimal reorder parameters at a given sampling time $t$, we simulate inventory within a forward-looking window from $t$ to $H$, where $H$ is the horizon. This simulation is done with the MRP system and needs to be repeated multiple times, incorporating different samples from the uncertainty distributions. We refer to this process as the forward-looking simulation, which includes two steps:
\begin{enumerate}
    \item Run a MRP simulation within the forward-looking window;
    \item Generate multiple uncertainty realizations of the simulated inventory by perturbing the demand forecast, the delay and quantity shortage of the arriving orders.
\end{enumerate}

The main inputs to the MRP simulation include:
\begin{enumerate}
    \item The demand forecast made on $t$ for the forward-looking window $C_{\textrm{DF}}^{t}\in\mathbb{R}^{H}$,
    \item The existing orders scheduled to arrive within a seeding window
    \begin{equation}
        A_{\textrm{seed}}^{t_s} = [a_{\textrm{seed}}^{t}, \dots, a_{\textrm{seed}}^{t+t_{s}-1}]\in\mathbb{R}^{t_s},
    \end{equation}
    where $t_{s}$ is the length of the seeding window, and $a_{\textrm{seed}}^{\tau}, \forall \tau\in[t, \dots, t+t_s-1]$ represents the seeded arrival on timestamp $\tau$,

    \item The candidate SSV and ST,
    \item The initial inventory position at sampling time $x_{0}^{t}$ (see Section.~\ref{subsec:detailed_model_treatements} for more details).
\end{enumerate}
Note that the size of the seeding window $t_s$ varies depending on the customers business logic, but oftentimes it is set to lead time ($LT$).

Given the new orders scheduled by MRP, we generate a total of $N_r$ realizations of inventory curves. The inventory at time $\tau+1$ for a realization $r$ is given by

\begin{equation}
\begin{aligned}
    x^{\tau+1}(r) & = x^{\tau}(r) - c_\text{A}^{\tau}(r) - c_\text{MM}^{\tau}(r) + s_s^{\tau}(r) \in\mathbb{R}^{1}\\
    & \forall \tau\in[t, \dots, t+H-1],
\end{aligned}
\end{equation}
where $x^{\tau}(r)$ is the inventory at time $\tau$ for the $r$-th realization, $c_\text{A}^{\tau}(r)$ is the sampled consumption at time $\tau$, $c_\text{MM}^{\tau}(r)$ is the sampled material movement at time $\tau$, $s_s^{\tau}(r)$ is sampled purchase order arrivals at $\tau$. Note that the inventory level $x$ here is no longer stochastic (not bold) with specific samples.

Specifically, the sampled consumption at time $\tau$ is given by
\begin{equation}
    c_\text{A}^{\tau}(r) = c_\textrm{DF}^{t, \tau} - u_\textrm{DF}(r)\in\mathbb{R}^{1},
\end{equation}
where $c_\text{DF}^{t, \tau}$ is the demand forecast made on time $t$ for time $\tau$, $u_\text{DF}(r)$ is the $r$-th sample from $U_\text{DF}$ (Eq.~\ref{equ:dfu}).

The sampled material movement at time $\tau$ is given by
\begin{equation}
    c_\text{MM}^{\tau}(r) = u_\text{MM}(r)\in\mathbb{R}^{1},
\end{equation}
where $u_\text{MM}(r)$ is the $r$-th sample from $U_\text{MM}$ (Eq.~\ref{equ:mm-uncertainty}).

Let ${o_1, o_2, ..., o_n}$ be the set of orders MRP placed,  $t_e^{o_j}$ and $s_e^{o_j}, \forall j\in[1,...,n]$ be the expected arrival time and expected quantity for order $o_j$. The $r$-th sampled arrival time for order $o_j$, denoted by $t_s^{o_j}(r)$ is given by
\begin{equation}
    t_s^{o_j}(r) = t_e^{o_j} + u_\text{ST}(r)\in\mathbb{R}^{1},
\end{equation}
where $u_\text{ST}(r)$ is $r$-th sample from $U_\text{ST}$ (Eq.~\ref{equ:squ}). The $r$-th sampled quantity for order $o_j$, denoted by $s_s^{o_j}(r)$ is given by
\begin{equation}
    s_s^{o_j}(r) = s_e^{o_j} + u_\text{SQ}(r)\in\mathbb{R}^{1},
\end{equation}
where $u_\text{SQ}(r)$ is the $r$-th sample from $U_\text{SQ}$ (Eq.~\ref{equ:stu}). Then the sampled purchase order arrivals at $\tau$, denoted by $s_s^{\tau}(r)$, is given by
\begin{equation}
    s_s^{\tau}(r) = \sum_{o_j} s_s^{o_j}(r)\in\mathbb{R}^{1}, \quad \forall t_s^{o_j}(r) = \tau.
\end{equation}

The inventory realizations generated by the forward-looking simulation can be expressed by
\begin{equation}
\label{equ:inventory_realizations}
\begin{aligned}
    X^{:H}(r) & = \left[x^{t}(r), \dots, x^{t+H-1}(r)\right] \in \mathbb{R}^{H}, \\
    & \forall r \in [1, \dots, N_r].
\end{aligned}
\end{equation}
Those realizations $X^{:H}(r), \forall r \in [1, \dots, N_r]$ are used to represent the stochastic variable $\mathbf{X}^{:H}$ and will be consumed by the k-iteration algorithm to produce recommendations for the reorder parameters.

\subsubsection{The k-iteration algorithm}
\label{subsec:k_iteration}
The k-iteration algorithm (also included in the reorder parameter estimation step shown in Fig.~\ref{fig:algorithm_workflow}) is used to map the hyper-parameters to the reorder parameter recommendations at each sampling time $t$. This algorithm is triggered in both training and validation (live production) phases of the algorithm workflow.

\begin{algorithm}
\caption{The k-iteration algorithm (for the Safety Stock MRP)}
\label{algo:k_ieration}
\begin{algorithmic}[1]
\State Set ST to the value corresponding to the STP (which is invariant of the k-th iteration)\;
\State Set $SSV = 0$, $h>0$\;
\While{k $<$ max iteration AND $h>0$}
    \State Get $N_r$ samples from each of the empirical uncertainty distribution (i.e., $U_{ST}$, $U_{SQ}$, $U_{MM}$, $U_{DF}$)\;
    \State Apply the forward-looking simulation to get $X^{:H}(r), \quad \forall r \in [1, \dots, N_r]$\;
    \State Calculate the inventory cost and the SL associated with each of the inventory realization $r$\;
    \If {Eq.~\ref{equ:k_iteration_stop_cond} does not hold}
        \State Compute the inventory deficit $h$ (Fig.~\ref{fig:k_iteration_compute_h})
    \Else
        \State Set $h=0$
    \EndIf
    \State Update $\textrm{SSV}_{k+1} = \textrm{SSV}_{k} + h$ and add incremental counter $k=k+1$\;
\EndWhile
\State \Return{(SSV, ST)}
\end{algorithmic}
\end{algorithm}

\begin{figure}[!htp]
  \centering
  \includegraphics[width=\linewidth]{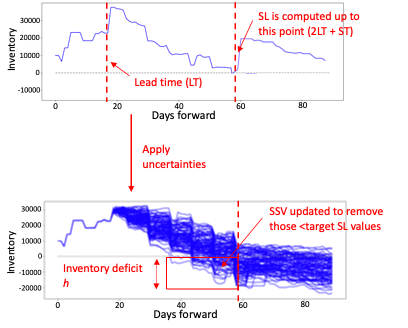}
  \caption{The k-iteration algorithm: compute inventory deficit $h$}
  \label{fig:k_iteration_compute_h}
\end{figure}

The procedure of the k-iteration algorithm can be found in Algorithm~\ref{algo:k_ieration}. At a sampling time $t$, the ST is determined by the STP and the empirical distribution of the supplier time uncertainty $U_{ST}$ in the corresponding USW. For SSV, it is initialized to 0. Then we enter the loop (i.e., the k-iteration) for increasing SSV incrementally till the $\textrm{SL}_{\textrm{min}}$ constraint in Eq.~\ref{equ:mpc_formulation} is met. In every iteration of this loop, a total of $N_r$ realizations of the simulated inventory curves are generated (i.e., $X^{:H}(r), \forall r \in [1, \dots, N_r]$) with the forward-looking simulation procedure. The inventory cost and SL associated with each of the realization of the simulated inventory (i.e., $X^{:H}(r)$) is computed. An inventory deficit $h$ is calculated to quantify how much lift of the inventory is needed to make sure $\psi(\mathbf{X}^{:H})$ is greater or equal to $\textrm{SL}_{\textrm{min}}$. For example, in Fig.~\ref{fig:k_iteration_compute_h}, if the minimum target SL and the SLP are both assumed to be 100\%, the inventory deficit $h$ would be the height of the red-boxed area, since that is how much the inventory level needs to be lifted so that 100\% (SLP) of the realization curves would satisfy a minimum target SL of 100\%. Note that for the MRP simulations in the k-iteration algorithm, the arrivals in $[0, LT)$ ($LT$ represents lead time) is seeded with the simulated arrivals from the MRP orders placed during the previous sampling time, since the order placed at the simulation of the current sampling time will not arrive the till after $LT$. Therefore, here the service level is computed in the time period of $[LT, \textrm{min}(2LT+ST, H)]$ due to the fact that the service level in $[0, LT)$ is not in any way affected by the choice of SSV at sampling time $t$, and the service level for later than $2LT+ST$ is likely overwritten in the next sampling time.

\begin{figure}[!htp]
  \centering
  \includegraphics[width=\linewidth]{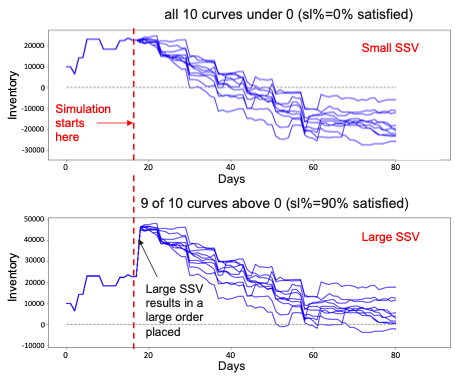}
  \caption{The k-iteration algorithm: lift SSV}
  \label{fig:k_iteration_lift_ssv}
\end{figure}

Once the inventory deficit $h$ is determined, the SSV for the next iteration would be calculated by adding $h$ to the current $\textrm{SSV}_{k}$ ($k$ denotes the iteration). The new $\textrm{SSV}_{k+1}$ would be used to run the MRP simulation and the uncertainty realizations for the next iteration $k+1$. As is shown in Fig.~\ref{fig:k_iteration_lift_ssv}, the lifted SSV is usually identical to the lift of inventory level since it would result in a large order arrival right after the lead time (LT) to lift the over all inventory level. The SLP will increase as a result. The iteration would continue till the max number of iteration is arrived or there is no need to lift anymore ($h=0$), which at each sampling time $t$, is equivalent to the following expression (Eq.~\ref{equ:k_iteration_stop_cond}) utilizing the notations developed in Eq.~\ref{equ:mpc_formulation} and \ref{equ:service_level}.

\begin{equation}
\label{equ:k_iteration_stop_cond}
    \gamma_{p_{\textrm{SL}}}\left[\frac{\sum_{i=\textrm{LT}}^{\textrm{min}(2\textrm{LT}+\textrm{ST},H)}\delta(\mathbf{x}^{t+i})}{\textrm{min}(\textrm{LT}+\textrm{ST}, H-\textrm{LT})}\right] - \textrm{SL}_{\textrm{min}} \geq 0,
\end{equation}
where $p_{\textrm{SL}}$ is the SLP.

Note that in the k-iteration procedure, we have configurable choice of the initial inventory, which we will discuss in more details in the next section.

\subsubsection{Other algorithm treatments}
\label{subsec:detailed_model_treatements}

There are a few key treatments of the workflow that would have significant impact on the performance and outputs of the algorithm. They include the batched uncertainty realizations with the MRP simulation, the ``steady-state" initial inventory of the k-iteration algorithm, the optimization frequency, and the demand forecast pre-processing. We will explain those treatments in this section.

\paragraph{Batched uncertainty realizations with MRP}\mbox{}\\

In the k-iteration algorithm described above, each uncertainty realization $r$ does not trigger another MRP simulation, which means for $N_r$ simulated inventory curves, the MRP simulation is only performed once. The rest of the curves are generated by simply perturbing the MRP result with the uncertainties sampled from the empirical distributions as a post-processing step. This process is referred to as the batched uncertainty realizations with MRP simulation. The benefit of such an approach is the reduction of computation time, since the uncertainty sampling is done much quicker in comparison to the MRP simulation. However, those realizations are referred to as the proxies since there could be difference between the results of this approach and those results if we sample uncertainty distributions during the MRP simulation processes (i.e., sequential uncertainty realizations with MRP).

The proxies of the simulated inventory is considered to be sufficiently accurate for the realizations of the simulated inventories, since those realizations are only used to determine how much SSV should be lifted, rather than accurately computing inventory cost and SL.

\paragraph{``Steady-state" initial inventory}\mbox{}\\

Intuitively, the initial inventory used for the k-iteration algorithm at sampling time $t$ is taken from the simulated inventory computed from the previous sampling time, shown in Eq.~\ref{equ:initial_inventory_simulated}. Note that at the beginning of the simulation period (for either training or validation), the simulated inventory is seeded with the actual inventory.

\begin{equation}
\label{equ:initial_inventory_simulated}
    x_{0}^{t} = x_{\text{sim}}^{t},
\end{equation}
where $x_{0}^{t}$ represents the initial inventory for the k-iteration algorithm starting at sampling time $t$. The variable $x_{\text{sim}}^{t}$ represents the inventory simulated at the previous sampling time.

The initial inventory determined by Eq.~\ref{equ:initial_inventory_simulated} is sensitive to the simulation process. Specifically, if the simulated inventory is high at $t$, the initial inventory of the k-iteration would be high. Therefore, all the realizations in the forward-looking simulation would be far above zero (satisfying the $\textrm{SL}_{\textrm{min}}$ constraint) even when we set SSV to 0. As a result, the algorithm will give a zero SSV recommendation as long as the current simulated inventory is high.

Ideally, the customers would expect the SSV recommendations be independent of the current inventory level. To achieve this goal, we propose an approach that always use a ``steady-state" inventory as the initial inventory fro the k-iteration algorithm. The initial inventory of the k-th iteration is denoted in Eq.~\ref{equ:initial_inventory_equation}.

\begin{equation}\label{equ:initial_inventory_equation}
    x^{t}_{0, k} = \text{SSV}_{k} + \sum^{\text{LT-1}}_{i=0}(c_{\text{DF}}^{t, t+i} - c_{\text{A}}^{t+i}),
\end{equation}
where $x^{t}_{0,k}$ represents the initial inventory for the $k$-th iteration starting at sampling time $t$. The variable $\text{SSV}_{k}$ is the SSV candidate at the $k$-th iteration. The variable $\text{LT}$ denotes the lead time. The variable $c_{\text{DF}}^{t, t+i}\in\mathbb{R}^{1}$ denotes the demand forecast for timestamp $t+i$ made at $t$, and $c_{\text{A}}^{t+i}\in\mathbb{R}^{1}$ represents the simulated arrivals at timestep $t+i$.

There are two notes for computing the initial inventory with Eq.~\ref{equ:initial_inventory_equation}. First, since the arrivals in $[0, LT)$ is seeded with the simulated arrivals, the inventory level at LT in the k-iteration procedure will be exactly at the level of the candidate SSV (i.e., $\text{SSV}_{k}$) regardless of the simulated inventory level at the sampling time $t$. Second, Eq.~\ref{equ:initial_inventory_equation} does not enforce the initial inventory to be larger or equal to zero, which means negative inventories could occur in $[0, LT)$ due to the choice of this initial inventory. However, since the inventory level will be brought back to positive at LT and the service level is only computed in $[LT, 2LT+ST]$, the potential negative inventories in $[0, LT)$ will not affect the SSV recommendations.

\paragraph{Optimization frequency}\mbox{}\\
Another key treatment is the frequency for which the reorder parameters are recommended. This frequency also corresponds to the frequency for which Eq.~\ref{equ:mpc_formulation} is solved. Ideally, solving Eq.~\ref{equ:mpc_formulation} at every sampling time $t$ would give an optimal result since new information is usually provided on a daily basis. However, it would be challenging (if possible at all) for the customers to adopt reorder parameter recommendations on a daily basis. Therefore, we designed a feature in the algorithm which allows the optimization frequency to be configurable by the user.

\begin{figure}[!htp]
  \centering
  \includegraphics[width=\linewidth]{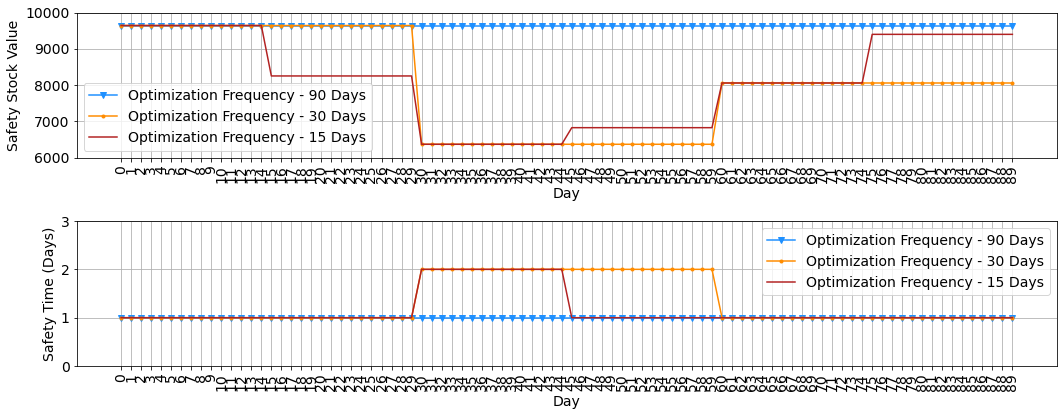}
  \caption{The SSV and ST recommendations with different optimization frequencies}
  \label{fig:ssv_st_opt_freq}
\end{figure}

If the optimization frequency in days is set to be a number $j$ larger than 1 day (e.g., 30 days). The reorder parameter recommendations are updated every $j$ days. For the days that those reorder parameters are not updated, they will default to the value computed at the previous sampling time. However, the MRP simulations in the order simulation step are still performed on a daily basis. Since the optimization frequency adds another constraint on the optimization of the reorder parameters, we state here that the more frequent the reorder parameters are allowed to change, the more room there will be for the algorithm to improve the results. To rationalize this statement, we setup experimental simulations with various optimization frequencies for a total length of 90 days. Fig.~\ref{fig:ssv_st_opt_freq} shows the SSV and ST recommendations with optimization frequencies of 90 days, 30 days, and 15 days. The blue lines are the SSV and ST recommendations with a frequency of 90 days, which remains unchanged during the entire simulation period. The yellow and red lines represents the reorder parameters with frequency of 30 and 15 days, respectively. The changes of lines in those colors are clearly observed.

\begin{figure}[!htp]
  \centering
  \includegraphics[width=\linewidth]{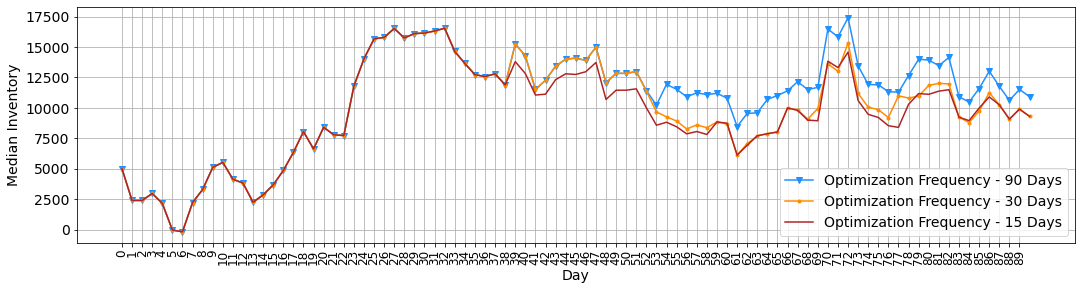}
  \caption{Simulated inventory levels with different optimization frequencies}
  \label{fig:inventory_opt_freq}
\end{figure}

Fig.~\ref{fig:inventory_opt_freq} demonstrates the median of the simulated inventory levels that correspond to each optimization frequency. In this example, the simulation inventory level corresponds to 90 days is the highest amount all three curves, and the one corresponds to 15 days is the lowest, while all three curves are 100\% above zero. The statement we made holds true in general, since the more frequent the reorder parameters are allowed to change, the less constraint there will be for the optimization, and therefore, a better result we would expect. To summarize, there is a clear trade-off between the algorithm performance and the customer adoption for the choice of optimization frequency. In reality, the decision made on this parameter would vary case-by-case.

\paragraph{Demand forecast pre-processing}\mbox{}\\
One more treatment to be discussed in this section is the pre-processing of the demand forecast data. The demand forecast provided by the customer usually comes with outliers. Those outliers are either due to input errors or inaccurate forecasts. The outliers usually are with value that is orders of magnitude larger than the regular data points, which could have large impact on the optimization result. To mitigate the impact of those outliers, we apply a clipping step on the demand forecast data, which is demonstrated in Eq.~\ref{equ:df_clipping}.

\begin{equation}
\label{equ:df_clipping}
    \tilde{C}_{\textrm{DF, 1D}} =
    \begin{cases}
        \eta_{C_{\textrm{DF}}}+n_c\cdot\sigma_{C_{\textrm{DF}}}, \\
         \quad \textrm{if} \quad C_{\textrm{DF, 1D}} > \eta_{C_{\textrm{DF}}}+n_c\cdot\sigma_{C_{\textrm{DF}}}, \\
        C_{\textrm{DF, 1D}}, \quad \textrm{otherwise}
    \end{cases}
\end{equation}
where $\tilde{C}_{\textrm{DF, 1D}}$ denotes the clipped demand forecast, and $C_{\textrm{DF, 1D}}$ is the raw demand forecast data computed in Eq.~\ref{equ:demand_forecast_1d}. The variables $\eta_{C_{\textrm{DF}}}$ and $\sigma_{C_{\textrm{DF}}}$ represent the median and the standard deviation of the demand forecast data, respectively. The variable $n_c$ is the a multiplier that controls the strictness of the clipping. It is a parameter that is configurable by the user. In our use cases shown in Section~\ref{section:results}, we default it to 5.0. The $\tilde{C}_{\textrm{DF, 1D}}$ is then used to calculate the convolved demand forecast $C_{\textrm{DF, 1D}}^{*}$ mentioned in Eq.~\ref{equ:demand_forecast_convolve}.

The outliers in the raw demand forecast data will also result in the over-estimation of the demand forecast uncertainty. Therefore, the demand forecast uncertainty derived from the clipped demand forecast will go through another step of clipping to further reduce those impacts. The clipping applied is shown in Eq.~\ref{equ:dfu_clipping}.

\begin{equation}
\label{equ:dfu_clipping}
    \tilde{U}_{\textrm{DF}} =
    \begin{cases}
        \eta_{U_{\textrm{DF}}}+n_u\cdot\sigma_{U_{\textrm{DF}}}, & \textrm{if} \quad U_{\textrm{DF}} > \eta_{U_{\textrm{DF}}}+n_u\cdot\sigma_{U_{\textrm{DF}}}\\
        U_{\textrm{DF}}, & \textrm{otherwise}
    \end{cases}
\end{equation}
where $U_{\textrm{DF}}$ is the demand forecast uncertainty defined in Eq.~\ref{equ:dfu} with $C_{\textrm{DF}}^{*}$ derived from $\tilde{C}_{\textrm{DF}}$ (instead of from $C_{\textrm{DF}}$ directly). The variables $\eta_{U_{\textrm{DF}}}$ and $\sigma_{U_{\textrm{DF}}}$ represent the median and the standard deviation of the demand forecast uncertainty, respectively. Similar to that for $\tilde{C}_{\textrm{DF}}$, the variable $n_u$ is the a multiplier that controls the strictness of the clipping. In our use cases shown in Section~\ref{section:results}, we default it to 1.0. Note that the movement uncertainty $U_{\textrm{MM}}$ can be treated in the same manner to generate clipped movement uncertainty $\tilde{U}_{\textrm{MM}}$.

\section{Inventory optimization at scale}
\label{section:platform}
\begin{figure*}[!htp]
  \centering
  \includegraphics[width=\textwidth]{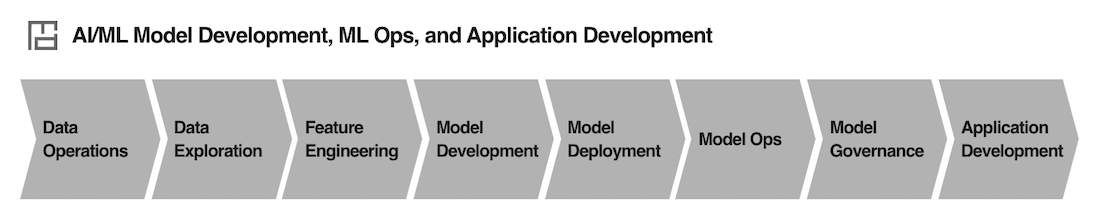}
  \caption{Flowchart for model deployment process in C3 AI Suite }
  \label{fig:model_deployment}
\end{figure*}

C3 AI Suite (the Platform), serves as the platform that hosts both data and computational frameworks that back the C3 Stochastic Inventory Optimization Application, is capable to provide enterprise users with data integration, machine learning development and operation tools, in a generalized application-development software stack (Figure~\ref{fig:model_deployment}).

Specific to the C3 SIO Application, we leveraged the Platform to integrate various formats of datasets from different data sources to the Suite such that the data can be easily consumed and orchestrate with the deployed optimization framework. The computational intensity from proposed framework can be easily identified based on the nature of multiple simulation runs. Typically, customers from different industries have hundreds of thousands to millions types of SKUs to optimize, leading to necessity of effective large scale optimization. The Model Deployment Framework in C3 AI Suite enable users to scale-out up to millions of ML models within a single application through pre-packaged automation including automatic scoring, compute resource auto-scaling, version management, alerting, hyper-parameter optimization, and auto retraining. By leveraging the model deployment framework on C3 AI Suite, it allows us to perform inventory optimization and simulation experiments across millions of SKUs, at scale and in a distributed fashion.

\subsection{Data integration}
\begin{figure*}[!htp]
  \centering
  \includegraphics[width=\textwidth]{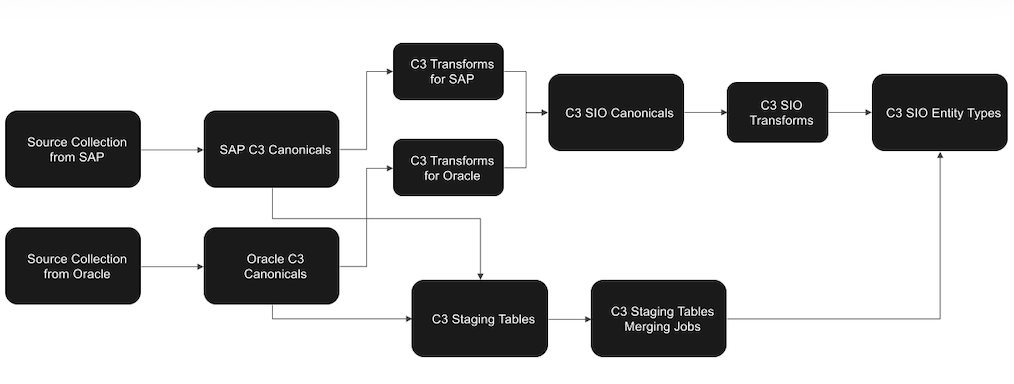}
  \caption{Flowchart for data integration process of ERP data}
  \label{fig:data_integration}
\end{figure*}

With the available data integration techniques supported on the Platform, we have configured the C3 Inventory Optimization Application to directly load the data from the enterprise resource planning (ERP) system to facilitate the data consumption of the SIO algorithm, as is shown in Figure~\ref{fig:data_integration}. The C3 Source Collection Types are used to map the source data locations and it can be configured to retrieve the source files on regular basis automatically with only retrieving the increments. The C3 Canonical Types are used for mapping the table headers of different data formats, and translating them into the data format that is native to the C3 AI Suite. Specifically to the C3 SIO Application, we have two layers of Canonicals, with the first layer mapping directly to the ERP source data, and the second layer mapping from a set that are transformed from the first layer, which follow closely to the business logic and are in closer form to the entity types in the Application. This provides users flexibility to provide data in either of the Canonical forms, depending on data availability. The C3 Transform Types are used for converting the data from Canonical Types and mapping the converted ones either to an additional layer of Canonicals or directly to entity types. Within the C3 Inventory Optimization Application, these entity types hold the information associated with certain business modules in resource planning, e.g. inventory level, material movement transactions, purchase orders, production orders, demand forecasts, etc.

\subsection{Model deployment}
By leveraging C3 AI Suite, we first incorporated various features based on the feature engineering tools and developed C3 AI ML pipeline to incorporate the whole simulation-optimization work flow. In order to deploy the ML pipelines on large number of SKUs, we created ML segments to support the model to SKUs relationships and then run distributed jobs to deploy ML pipelines on a large number of SKUs or a group of SKUs that belong to one segment. We can also configure champion and challenger model to make sure the live model is always providing the optimal results. Model operations and model governance are critical steps in managing the model recommendation and explaining the model predictions.

\section{Case study results}
\label{section:results}

\nomenclature[A]{SKU}{Stock keeping unit}

The C3 Stochastic Inventory Optimization algorithm has been tested on customer data from a variety of verticals. Here, we show the performance of C3 SIO with different levels of details and demonstrate the inventory savings that can be achieved by applying our algorithm.

In Section~\ref{subsec:results_sampled_SKU}, we talk about two representative SKUs come from one of our clients in the energy and manufacturing industry. We demonstrate in detail how the algorithm manages to recommend a lower SSV while keeping a satisfying service level. In the Section~\ref{subsec:results_sensitivity}, we expand the dataset and use 54 representative SKUs selected by the same client. We discuss the impact of some of the model parameters on the performance of the algorithm. In Section~\ref{subsec:results_verticals}, we illustrate the generality of the algorithm by showing its performance on the customer data across a variety of verticals.

\subsection{Optimization results for sampled SKUs}
\label{subsec:results_sampled_SKU}

In this part, we apply C3 SIO on two representative SKUs (i.e., SKU A and SKU B) that come from one of our clients in the energy and manufacturing industry. The minimum target service level for SKU A and SKU B are set to 91\% and 96\%, respectively. We have the data available for both SKUs from 2018-06-01 to 2020-10-01. The training period for C3 SIO starts at 2018-06-01 and ends at 2020-03-01, while the validation period starts at 2020-04-01 and ends at 2020-10-01. Since 2020-10-01 is the last date with data available, the operation date (or, live production date) is set to 2020-10-01. The optimization frequency is set to 30 days, which means our algorithm will be able to recommend new reorder parameters every 30 days.

According to the requirements of the client, we did not change STP and always set $\textrm{ST}=0$. Therefore, for the grid search in the training period, we will only select the optimal SLP from the SLP candidate list. The SLP candidate list is set to a default value of $[50,70,90,92.5,95]$ here. The default multiplier for clipping the demand forecast ($n_c$ in Eq.~\ref{equ:df_clipping}) is set to 5.0. The default multiplier for clipping the demand forecast uncertainty ($n_u$ in Eq.~\ref{equ:dfu_clipping}) is set to 1.0. The default minimum length of uncertainty sampling window $L_{\textrm{min}}^{\textrm{USW}}$ (Eq.~\ref{equ:usw-length}) is 30 days and the USW buffer $b^{\textrm{USW}}$ (Eq.~\ref{equ:usw-length}) is set to 14 days. We set uncertainty realizations in the forward-looking simulation procedure ($N_r$ in Eq.~\ref{equ:inventory_realizations}) to 100. We will alter those default values later in Section~\ref{subsec:results_sensitivity} to evaluate their impact on the performance of the algorithm.

Note that during the order simulation step in the C3 SIO workflow, we also need to generate multiple realizations (we note it as $N_{os}$ here) of the inventory to compare with the actual inventory the customer has. The variable $N_{os}$ is set to 10 in this case. Also, the client intends to focus on the savings for the inventory holding cost. Therefore, the variable $\pi_{o,t}$ in Eq.~\ref{equ:mpc_formulation} is set to 0 in this case.

\nomenclature[V]{$N_{os}$}{Number of uncertainty realizations in the order simulation step}

\begin{figure}[!htp]
  \centering
  \includegraphics[width=\linewidth]{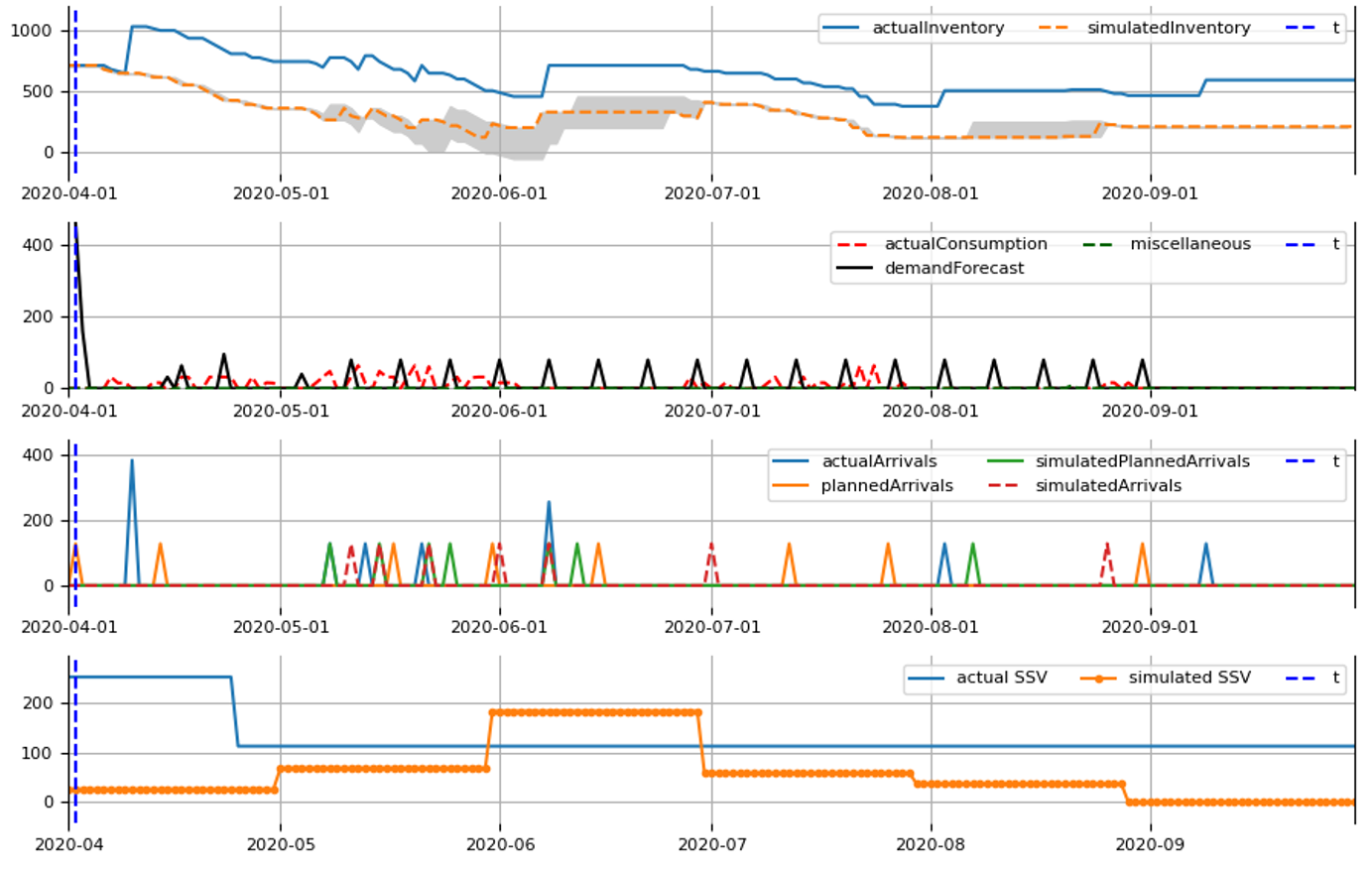}
  \caption{Simulation result for SKU A with Plot 1-4 from top to bottom. Plot 1: The simulated inventory and the actual inventory. Plot 2: The demand forecast and the actual consumption. Plot 3: Planned arrivals and the actual arrivals for both simulated and the actual inventories. Plot 4: The simulated SSV recommendations and the actual SSV.}
  \label{fig:result_fig_1}
\end{figure}

Figure~\ref{fig:result_fig_1} shows the results for the inventory optimization on SKU A. Plot 1 on the top of Fig.~\ref{fig:result_fig_1} shows the comparison between the median of the simulated inventory realizations (dashed orange curve) and the actual inventory (blue curve) over the validation period. The gray area surrounding the orange curve is the uncertainty range generated from the 10 uncertainty realizations. Plot 2 shows the demand forecast (black curve) and the actual consumption (red dashed curve). Plot 3 shows comparisons for the planned arrivals and the actual arrivals in the simulation and in actual, respectively. Plot 4 of Fig.~\ref{fig:result_fig_1} compares the simulated SSV recommendations (orange curve) with the actual SSV (blue curve).

Looking at Plot 1 in Fig.~\ref{fig:result_fig_1}, it is obvious that the simulated inventory in the dashed orange curve is overall at a lower level in comparison with the actual inventory (blue curve). The gray area shows the uncertainty range for the simulated inventory with 10 uncertainty realizations. As we can see, a small portion of the gray area drops below zero around June 1st, 2020, indicating a slight out-of-stock risk (potentially caused by some outliers in the uncertainties). However, the dashed orange curve (the median of all the curves in the gray area) is far away from 0, which gives us more confidence in the overall healthiness of the inventory level.

The black curve and the dashed red curve in the Plot 2 of Fig.~\ref{fig:result_fig_1} indicates that the demand forecast is significantly more than the actual consumption. We refer to this scenario as the over-forecast. Note that Plot 2 only shows the demand forecast made on April 1st, 2020. If we look into the demand forecasts on other days in the validation period, we will find that they are consistently over forecasting.

Plot 3 in Fig.~\ref{fig:result_fig_1} represents the arrival information for both simulated inventory and the actual inventory. The actual arrival is shown in blue, which corresponds to the increases of the blue curve in Plot 1. Comparing the blue curves and the orange curves (planned arrivals) in Plot 3, we will have a sense of the supplier time uncertainty (how much the actual arrivals are late than the planned arrivals) and the supplier quantity uncertainty (how many quantities the actual arrivals are less than the planned arrivals) for the actual inventory. Similarly, the dashed red curve in plot 3 corresponds to the simulated arrivals, which corresponds to the increases of the dashed orange curve in Plot 1. Comparing the red dashed curve and the green curve, we will know about the supplier time and quantity uncertainties for the simulated inventory.

Plot 4 in Fig.~\ref{fig:result_fig_1} shows the simulated SSV recommendation (orange) and the actual SSV (blue). The simulated SSV is overall much lower than the actual SSV resulting in less order placed in the validation period. This is the main reason why the overall simulated inventory level is lower than the actual in Plot 1. Note that the orange curve is stair-like because we set the optimization frequency to 30 days.

\begin{figure}[!htp]
  \centering
  \includegraphics[width=\linewidth]{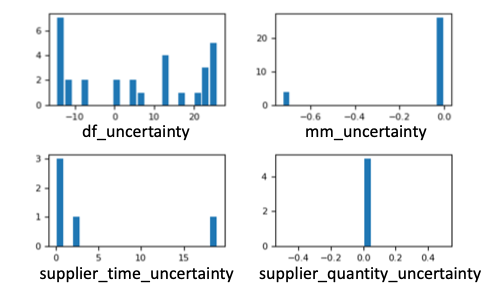}
  \caption{Sampled uncertainties within the USW ending at 2020-04-01 for SKU A. From upper left to lower right: demand forecast uncertainty, movement uncertainty (which includes both the miscellaneous and blocked movements), supplier time uncertainty and the supplier quantity uncertainty.}
  \label{fig:uncertainty_plot_1}
\end{figure}

To find the reason behind the low SSV recommendation shown in Plot 4 of Fig.~\ref{fig:result_fig_1}, we need to look into the uncertainty sampling process and the k-iteration algorithm. Fig.~\ref{fig:uncertainty_plot_1} demonstrates the uncertainties sampled at the beginning of the validation period. Note that the lead time for this SKU is 6 days, with the USW buffer $b^{\textrm{USW}}$ set to 14 days, the uncertainty sampling window spams from 2020-03-12 to 2020-04-01. The plot on the left in Fig.~\ref{fig:result_fig_1} shows the histogram of the demand forecast uncertainty. There are more positive than the negative samples in the histogram, which is consistent with the over-forecast we observed in Plot 2 of Fig.~\ref{fig:result_fig_1}.

\begin{figure}[!htp]
  \centering
  \includegraphics[width=\linewidth]{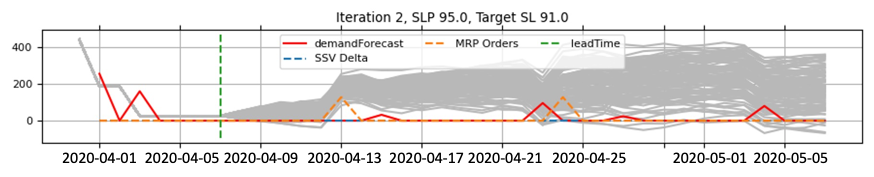}
  \caption{The forward-looking simulation result in the first iteration of the k-iteration procedure at 2020-04-01 for SKU A}
  \label{fig:k_iteration_result_1}
\end{figure}

In the k-iteration algorithm, the demand forecast is used as the input for the MRP run. Then, the output of the MRP run is ``corrected'' to the actual consumption by sampling the demand forecast uncertainty. If we have more positive samples in the demand forecast uncertainty, it means the actual consumption is more likely to be less than the demand forecast used as the input of the MRP. Therefore, the simulated inventory level is more likely to be pushed up in the uncertainty realizations. Fig.~\ref{fig:k_iteration_result_1} shows the simulation with uncertainty realizations in the first iteration of the k-iteration procedure. The gray curves represent the uncertainty realizations of the simulated inventory. As we can see, most of the uncertainty realizations are above 0 due to the reason we discussed above. Therefore, the inventory deficit is small and we do not need to lift too much on SSV to satisfy the constraint on the service level.

As we can see in this example, C3 SIO algorithm offered an objective view for potential risks involved every time when we recommend the SSV. This view is obtained from MRP run with multiple uncertainty realizations in the k-iteration procedure. With SKU A, such a view helped us reduced the SSV significantly compared to the actual SSV used by the customer.

\begin{figure}[!htp]
  \centering
  \includegraphics[width=\linewidth]{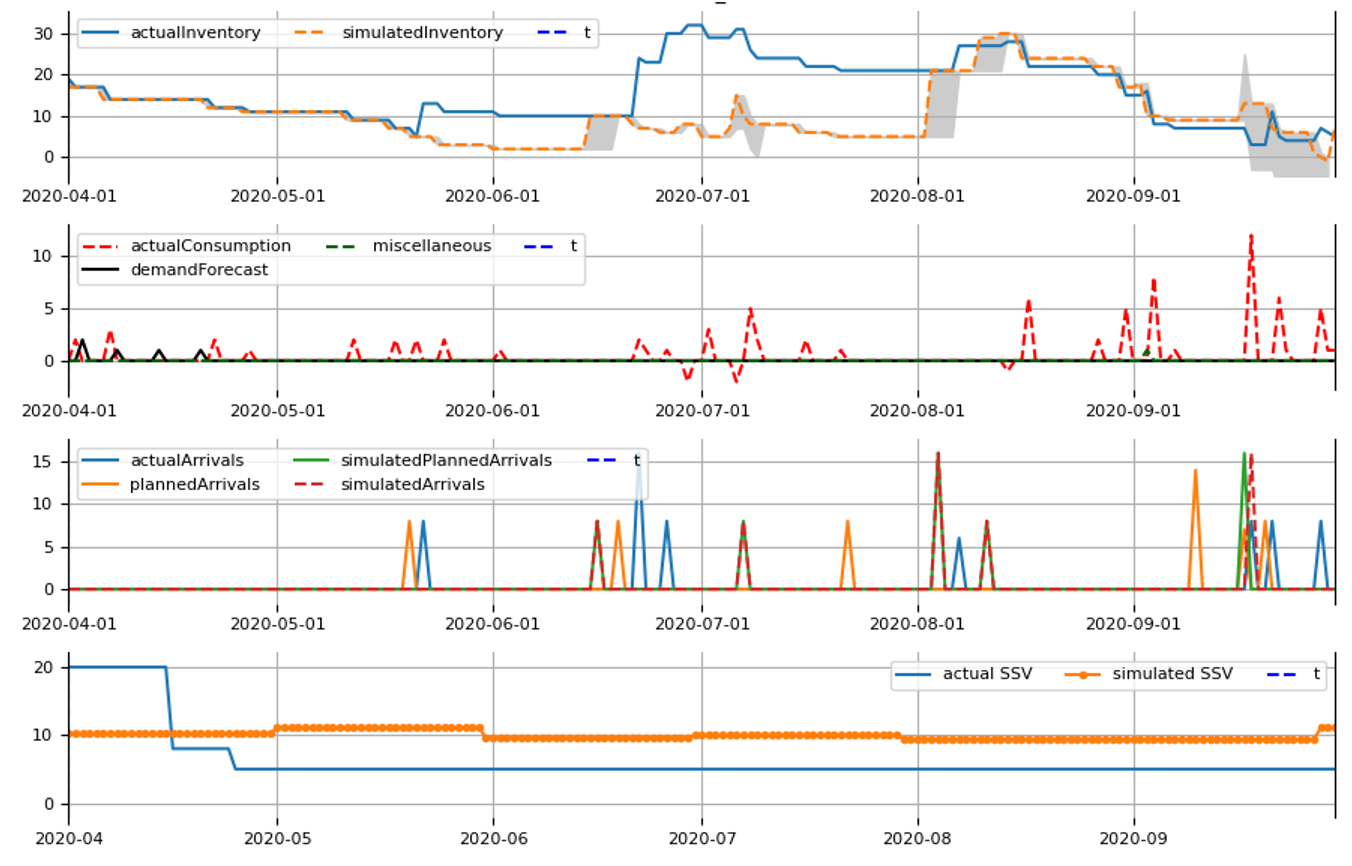}
  \caption{Simulation result for SKU B with Plot 1-4 from top to bottom. Plot 1: The simulated inventory and the actual inventory. Plot 2: The demand forecast and the actual consumption. Plot 3: Planned arrivals and the actual arrivals for both simulated and the actual inventories. Plot 4: The simulated SSV recommendations and the actual SSV.}
  \label{fig:result_fig_2}
\end{figure}

Now let's look at another SKU with different characteristics. Fig.~\ref{fig:result_fig_2} shows the result for the inventory optimization on SKU B, which is similar to Fig.~\ref{fig:result_fig_1} for SKU A with all the plots revealing the same information. We want to point out that there is one big difference for SKU B compared with SKU A regarding the demand forecast. As we can see in Plot 2 of Fig.~\ref{fig:result_fig_2}, the black curve is significantly less than the dashed red curve, which indicates a significant under-forecast for the demand forecast. Note that the black curve in Plot 2 only shows the demand forecast made on April 1st, 2020. In fact, the demand forecast made on most of the days in the validation period are all zero. This is not a rare case in the customer's data. We observe significant under-forecast or complete missing of demand forecasts for the vast majority of the in scope SKUs provided by the customer, due to the fact that the customer lacks the capability to predict the consumption for the SKUs it holds. Nevertheless, as we will demonstrate with SKU B here, C3 SIO algorithm is capable to maintaining the service level with SSV recommendations even when the demand forecast is very sparse or completely missing.

Plot 4 in Fig.~\ref{fig:result_fig_2} indicates that the algorithm overall recommends a much higher SSV compared to the actual SSV consumed by the customer. Yet, Plot 1 shows that that higher SSV recommendation ends up with a overall lower simulated inventory level compared to the actual inventory, while still satisfies the service level constraint. It is not surprising that the actual SSV is much lower than the simulated SSV, but the actual inventory ends up in a much higher position compared to the simulated inventory. Since the actual inventory is managed by the human operator, and they may not place order exactly following the MRP output. In the case of SKU B, it is likely that the human operator placed more orders than what MRP recommends. Should the orders be placed exactly following the actual SSV, the actual inventory will have significant service level issue (i.e., serious stock-out).

\begin{figure}[!htp]
  \centering
  \includegraphics[width=\linewidth]{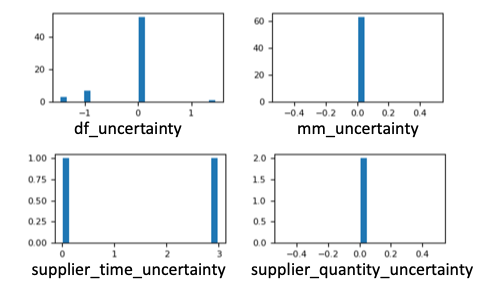}
  \caption{Sampled uncertainties within the USW ending at 2020-04-01 for SKU B. From upper left to lower right: demand forecast uncertainty, movement uncertainty (which includes both the miscellaneous and blocked movements), supplier time uncertainty and the supplier quantity uncertainty.}
  \label{fig:uncertainty_plot_2}
\end{figure}

It is quite intuitive that when the demand forecast is missing, we maintain a high SSV level to account for the unexpected consumption. Fig.~\ref{fig:uncertainty_plot_2} shows the uncertainties sampled at the beginning of the validation period. The left plot with the demand forecast uncertainty indicates an overall under-forecast, since most of the non-zeros samples are negative. The negative samples in the demand forecast uncertainty would pull the inventory level lower in the uncertainty realizations of the k-iteration procedure, forcing a higher SSV to maintain the service level.

\begin{figure}[!htp]
  \centering
  \includegraphics[width=\linewidth]{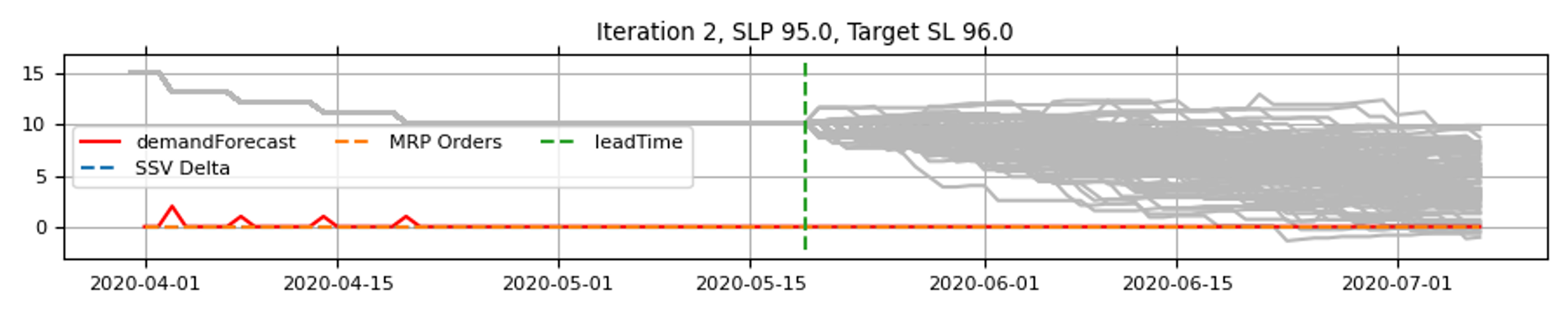}
  \caption{The forward-looking simulation result in the final iteration of the k-iteration procedure at 2020-04-01 for SKU B}
  \label{fig:k_iteration_result_2}
\end{figure}

Fig.~\ref{fig:k_iteration_result_2} shows the uncertainty realizations in the last iteration of the k-iteration procedure on 2020-04-01. The gray curves in the plot are the uncertainty realizations. Noticing that the overall downward trend of the inventory is due to the negative samples of the demand forecast uncertainty. It indicates that the demand forecast fed to the MRP run is an under-forecast, so the ``corrections" are applied to pull the inventory level down. Therefore, to meet the service level constraint, a higher SSV is recommended.

\subsection{Sensitivity study on the model parameters}
\label{subsec:results_sensitivity}

Different from the hyper-parameters, which specifically refer to SLP and STP. There are a few model parameters which are also critical to the performance of C3 SIO. Those model parameters include the multiplier for clipping the demand forecast (i.e., the variable $n_c$ in Eq.~\ref{equ:df_clipping}), the multiplier for clipping the demand forecast uncertainty (i.e., the variable $n_u$ in Eq.~\ref{equ:dfu_clipping}), the candidate list for SLP, the length of the offset (buffer) for USW (i.e., the variable $b^{\textrm{USW}}$ in Eq.~\ref{equ:usw-length}), and the number of uncertainty realizations (i.e., the variable $N_r$ in Eq.~\ref{equ:inventory_realizations}).

To evaluate the impact of the model parameters on the model performance, we expand  the experiment subjects to an ensemble of 54 SKUs, which comes from same client as the SKUs shown in Section~\ref{subsec:results_sampled_SKU}. Those 54 representative SKUs are hand-picked by the client's subject matter experts with an intention to cover various scenarios they encounter in the business operation.

For each of the study, we compare a few key performance metrics including the percentage of SKUs meeting the target service level (i.e., service level adherence rate, $r_\textrm{AD}$), inventory holding cost saving percentage averaged over the entire validation period ($\bar{s}_{\textrm{inv}}$), safety stock holding cost saving percentage averaged over the entire validation period ($\bar{s}_{\textrm{ss}}$), safety stock holding cost saving percentage at the operation date ($s_{\textrm{ss, op}}$).

The service level adherence rate is defined as follows:

\begin{equation}
\label{equ:service_level_adherence_rate}
    r_{\textrm{AD}} = \frac{\sum_{i=1}^{N_{\textrm{SKU}}}\delta(\textrm{SL}_{i} - \textrm{SL}_{i, \textrm{target}}) }{N_{\textrm{SKU}}},
\end{equation}
where variable $N_{\textrm{SKU}}$ represents the total number of SKUs in the scope. In our experiment setting, $N_{\textrm{SKU}} = 54$. The notation $\delta(\cdot)$ is a level set function defined in Eq.~\ref{equ:level_set}. The variables $\textrm{SL}_{i}$ and $\textrm{SL}_{i, \textrm{target}}$ denotes the service level from C3 SIO simulation and the target service level for SKU $i$, respectively.

The inventory holding cost saving percentage averaged over validation period is defined as follows:
\begin{equation}
\label{equ:inv_saving_val_mean}
    \bar{s}_{\textrm{inv}} = \frac{\sum_{t=1}^{N_{T}}(x_{a}^{t} - x_{\textrm{sim}}^{t})\pi_{h}^{t}}{\sum_{t=1}^{N_{T}}x_{a}^{t}\pi_{h}^{t}},
\end{equation}
where the variable $t$ denotes the timestep in the validation period, while $N_{T}$ represents the total number of days in the validation period. The variables $x_{a}^{t}$ and $x_{\textrm{sim}}^{t}$ depict the actual and the simulated inventory level at timestamp $t$, respectively, and $\pi_{h}$ represents the inventory holding cost.

Similarly, the safety stock holding cost saving percentage averaged over validation period, and the safety stock holding cost saving percentage at operation date are defined in Eqs.~\ref{equ:ssv_saving_val_mean} and \ref{equ:ssv_saving_val_op} below.
\begin{equation}
\label{equ:ssv_saving_val_mean}
    \bar{s}_{\textrm{ss}} = \frac{\sum_{t=1}^{N_{T}}(\textrm{SSV}_{a}^{t} - \textrm{SSV}_{\textrm{sim}}^{t})\pi_{\textrm{h}}^{t}}{\sum_{t=1}^{N_{T}}x_{a}^{t}\pi_{\textrm{h}}^{t}},
\end{equation}

\begin{equation}
\label{equ:ssv_saving_val_op}
    s_{\textrm{ss, op}} = \frac{N_{T}(\textrm{SSV}_{a}^{\textrm{op}} - \textrm{SSV}_{\textrm{sim}}^{\textrm{op}})\pi_{h}^{\textrm{op}}}{\sum_{t=1}^{N_{T}}x_{a}^{t}\pi_{h}^{t}},
\end{equation}
where the notation $\textrm{SSV}$ represents the safety stock value, and the superscript $\textrm{op}$ represents the timestep at the operation date, which usually is the date right after the validation period finishes. Note that in Eqs.~\ref{equ:ssv_saving_val_mean} and \ref{equ:ssv_saving_val_op}, the denominator is the total inventory holding cost (rather than the total safety stock holding cost) over the validation period. This is worth to note that $\bar{s}_{\textrm{ss}}$ and $s_{\textrm{ss, op}}$ are more important to this client compared to $\bar{s}_{\textrm{inv}}$ in this case.

\begin{figure}[!htp]
  \centering
  \includegraphics[width=\linewidth]{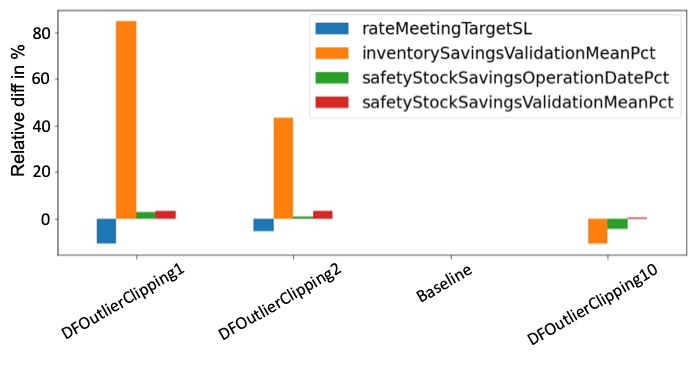}
  \caption{Sensitivity study on the clipping multiplier for raw demand forecast data}
  \label{fig:sa_DFOutlierClipping}
\end{figure}

The results of the sensitivity study for $n_c$ (the multiplier for clipping the demand forecast) is shown in Fig.~\ref{fig:sa_DFOutlierClipping}. The we set $n_c = 5.0$ as the baseline case, and experimented with both $n_c$ larger and smaller than the baseline. The x-axis shows $n_c$ of 1.0, 2.0, 5.0, and 10.0, successively. The blue bars represent $r_{\textrm{AD}}$, the orange bars represent $\bar{s}_{\textrm{inv}}$, the green bars represent $s_{\textrm{ss, op}}$, and the red bars represent $\bar{s}_{\textrm{ss}}$. With $n_c$ larger than the baseline, we see a clear dropping for all the saving metrics (i.e., $\bar{s}_{\textrm{inv}}$, $\bar{s}_{\textrm{ss}}$, $s_{\textrm{ss, op}}$) with no improvement at $r_\textrm{AD}$. On the other hand, if we make the clipping criteria more strict and have $n_c$ smaller than the baseline, we can observe a clear improvement of the saving metrics, especially for the inventory holding cost savings (orange bars). However, the service level adherence rate will be lower.

Obviously, there is a trade-off for between the service level and the saving. Starting from $n_c=1.0$, if we aim to improve $r_{\textrm{AD}}$ by increasing $n_c$, the improvement would come at the cost of reducing the savings. However, after $n_c = 5.0$, increasing $n_c$ will no longer provide improvement on $r_{\textrm{AD}}$. Therefore, with the presented data, we would stop at $n_c=5.0$ to avoid additional reduction on the savings.

\begin{figure}[!htp]
  \centering
  \includegraphics[width=\linewidth]{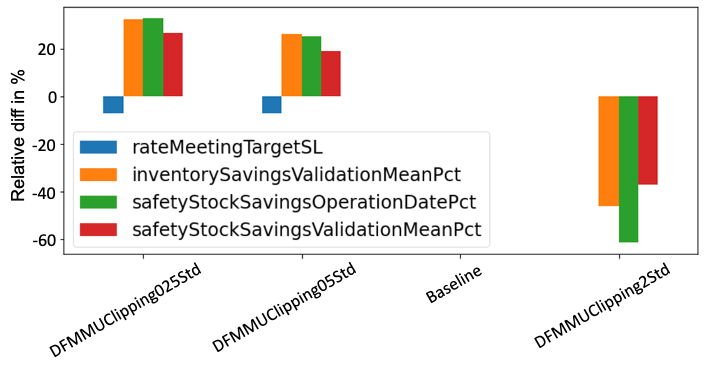}
  \caption{Sensitivity study on the clipping multiplier for demand forecast uncertainty and movement uncertainty}
  \label{fig:sa_DFMMUClipping}
\end{figure}

Similar study can be applied to $n_u$ (the multiplier for clipping the demand forecast uncertainty). Fig.~\ref{fig:sa_DFMMUClipping} compares the service level adherence rate and the saving metrics for inventory optimization results with various $n_u$. The baseline for $n_u$ is 1.0. The bars map to the same metrics as Fig.~\ref{fig:sa_DFOutlierClipping}. We compared cases with $n_u$ set to 0.25, 0.5 and 2.0. As we increase the $n_u$, the $r_{\textrm{AD}}$ will improve at the cost of reducing savings. As $n_u$ gets larger than 1.0, the improvement for $r_{\textrm{AD}}$ stopped. The results here indicates that $n_u = 1.0$ is a reasonable choice, if improving $r_{\textrm{AD}}$ is more important increasing savings.

\begin{figure}[!htp]
  \centering
  \includegraphics[width=\linewidth]{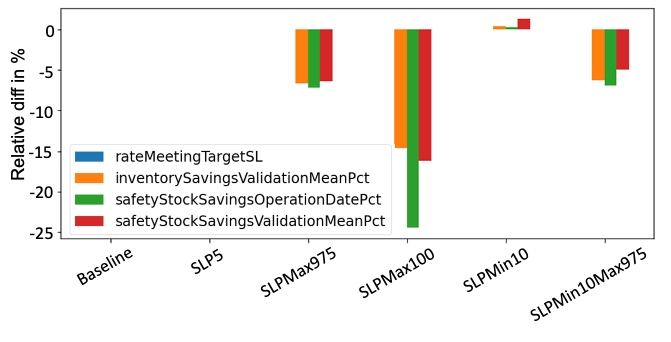}
  \caption{Sensitivity study on the candidate list for service level percentile (SLP)}
  \label{fig:sa_SLP}
\end{figure}

The hyper-parameter search results depends on the initial candidate list we provided to the algorithm. For example, if the grid search of SLP starts with different candidate SLP lists, the optimized SLP will be different and the resulting inventory or SSV savings would be affected. Fig.~\ref{fig:sa_SLP} demonstrates the results for the sensitivity study on the SLP candidate lists, since in this study we only optimizing on SLP as hyper-parameter. Each item on the x-axis notations a SLP candidate list. The mapping between the x-axis notation in Fig.~\ref{fig:sa_SLP} and the SLP candidate list can be found in Table~\ref{tab::slp_candidates}. For example, ``Baseline" corresponds to the initial SLP candidates of 50\%, 70\%, 90\%, 92.5\%, 95\%. For experiments with ``Baseline", the best SLP will be chosen from the above 5 candidates.

\begin{table}[!htb]
\centering
\begin{tabular}{l@{\hskip .1in} l}
    \hline
    Notation    &SLP candidates \\
        \\[-1em]
    \hline
    Baseline    &[50, 70, 90, 92.5, 95]\\
        \\[-1em]
    \hline
    SLP5        &[50, 55, 60, 65, 70, 75, 80,\\
                &85, 90, 92.5, 95]\\
        \\[-1em]
   \hline
    SLPMax975   &[50, 70, 90, 92.5, 95, 97.5]\\
        \\[-1em]
    \hline
    SLPMax100       &[50, 70, 90, 92.5, 95, 97.5, 100]\\
        \\[-1em]
    \hline

    SLPMin10        &[10, 30, 50, 70, 90, 92.5, 95]\\
        \\[-1em]
    \hline

    SLPMin10Max975  &[10, 30, 50, 70, 90, 92.5, 95, 97.5]\\
        \\[-1em]
    \hline
\end{tabular}
\caption{SLP candidates shown in Fig.~\ref{fig:sa_SLP}}
\label{tab::slp_candidates}
\end{table}

Figure~\ref{fig:sa_SLP} clearly indicates that the service level adherence rate is not affected as we change the SLP candidate list, for the 54 SKUs tested. In addition, as we increase the upper bound of the candidate list from 95\% to 97.5\% and 100\%, the saving metrics drop significantly. It is expected since a high SLP value means the recommended SSV will be more likely to be affected by the outliers of the uncertainties. For example, if we have an outlier in the supplier quantity uncertainty (a large shortage), setting SLP equal to 100\% means the recommended SSV will ensure the inventory meets target service level even when this outlier is present, which will end up with a very high SSV recommendation. On the other hand, if the lower bound of the candidate list is extended, we would see a slight increase of the savings. Since selecting a very low SLP candidate involves larger risks for service level hits, the decision of extending the lower bound or not will be made by the customers based on their risk profiles.

\begin{figure}[!htp]
  \centering
  \includegraphics[width=\linewidth]{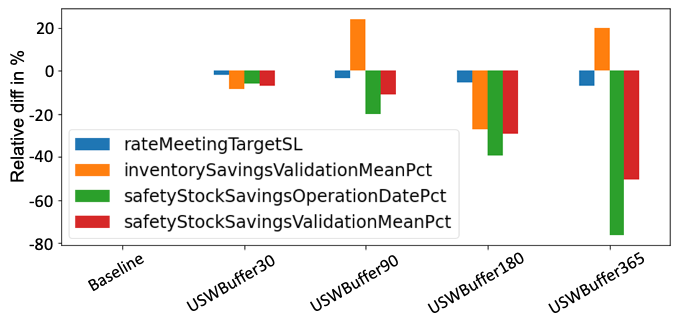}
  \caption{Sensitivity study on the buffer for the uncertainty sampling window (USW) size}
  \label{fig:sa_USWBuffer}
\end{figure}

Figure~\ref{fig:sa_USWBuffer} shows the results for cases with various buffers for the USW ($b^{\textrm{USW}}$). The baseline buffer is set to 14 days. Then the results with buffer of 30, 90, 180 and 365 days are plotted successively in the figure. Despite some minor inconsistencies between the inventory saving and the SSV saving, we observe a general downward trend for both saving metrics and the service level adherence rate $r_{\textrm{AD}}$. This result is expected since as the USW gets longer, outliers are more likely to appear in the sampling window, causing problems for both savings and service level. Therefore, it is a wise choice to keep a shorter USW as long as enough samples could be collected within the window.

\begin{figure}[!htp]
  \centering
  \includegraphics[width=\linewidth]{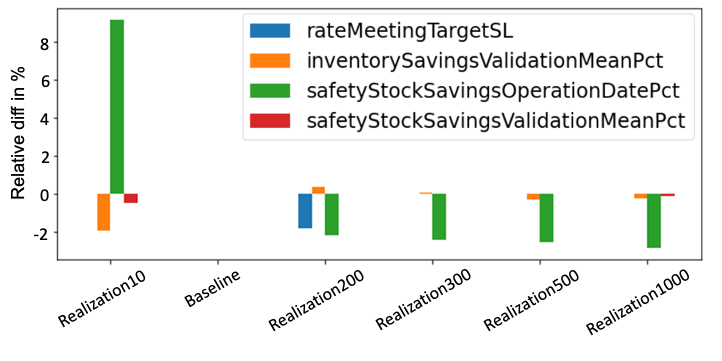}
  \caption{Convergence study on the number of uncertainty realizations in the forward-looking simulation procedure}
  \label{fig:sa_realizatioins}
\end{figure}

Figure~\ref{fig:sa_realizatioins} shows a convergence study on the number of uncertainty realizations in the forward-looking simulation ($N_r$). The baseline in the figure represents 100 realizations. As is observed from the figure, results experience large oscillation with smaller number of realizations. As the number of realizations approaches 500, both saving metrics and the service level adherence rate $r_{\textrm{AD}}$ stabilize.

\subsection{Summary of results for customers in different verticals}
\label{subsec:results_verticals}

In addition to the customer in the energy and manufacturing industry, we have applied C3 SIO algorithm to other customer data in a variety of verticals including building, manufacturing, healthcare, and energy in general. The results of the algorithm are summarized in Tab.~\ref{tab:verticals}. The columns in Tab.~\ref{tab:verticals} refers to the total number of SKUs in-scope, the actual inventory cost for the customer for the validation period, the simulated inventory cost generate by C3 SIO in the validation period, the actual SSV cost for the customer, the cost of SSV recommended by C3, the percentage of savings (either with inventory savings or SSV savings), the actual service level maintained by the customer, and the service level for the simulated inventory (with SSV recommended by C3).

Note that for some of the customers (i.e., the one in energy \& manufacturing vertical), they measure the savings by the percentage reduction on the SSV cost (defined in Eq.~\ref{equ:ssv_saving_val_mean}). We want to distinguish it from other inventory savings while interpreting the results.

We demonstrated in Tab.~\ref{tab:verticals} that C3 SIO algorithm is effective for large-scale inventory optimization across a variety of verticals. It is capable of reducing the inventory and the safety stock costs by 10-20\% while maintaining the service level similar to its current position.

\begin{table*}[!htb]
\centering
\begin{tabular}{l@{\hskip .1in} l@{\hskip .1in} l@{\hskip .1in} l@{\hskip .1in} l@{\hskip .1in} l@{\hskip .1in} l@{\hskip .1in} l@{\hskip .1in} l@{\hskip .1in} l}
    \hline
    \hline
    \textbf{Customers}   & \#SKUs  &$\text{Inv}_{a}$        &$\text{Inv}_{\text{C3}}$  &$\text{SSV}_{a}$        &$\text{SSV}_{\text{C3}}$  &Saving \%     &$\text{SL}_{a}$        &$
    \text{SL}_{\text{C3}}$ \\
        \\[-1em]
    \hline
    \hline
    Customer 1  &501    &\$10M      &\$7.5M    &-      &-   &25\% &99\%  &98.0\% \\
        \\[-1em]
   \hline
    Customer 2    &20.6k   &\$64M &- &\$11M &\$1.9M &14.2\%* &99.9\%** &98.7\%  \\
        \\[-1em]
    \hline
    Customer 3    &1.4k &\$24M &\$21M &- &- &14\% &88.6\% &90.3\%    \\
        \\[-1em]
    \hline

    Customer 4  &1.3K &\$7.7M       &\$6.5M     &-   &- &15\%   &97.6\% &91.0\%  \\
        \\[-1em]
    \hline
    Customer 5  &1k  &\$11.7M        &\$9.2M     &- &-   &21.2\% &98.6\% &94.2\% \\
        \\[-1em]
    \hline
    \hline

\end{tabular}
\caption{Inventory optimization results for different customers (k: thousand, M: million) with EUR converted to USD at a rate of 1.1. *SSV saving percentage is calculated following Eq.~\ref{equ:ssv_saving_val_mean}. **The actual SL is unreasonably high due to human intervention with extra information}
\label{tab:verticals}
\end{table*}

\section{Conclusion and discussions}
In this paper, we introduced the challenges for the large-scale supply chains, and proposed to solve them with C3 Stochastic Inventory Optimization (SIO) algorithm. The C3 SIO algorithm works well with a variety of MRP systems, and we demonstrated in detail how it works with a specific type of MRP system (i.e., the safety stock MRP). The C3 SIO algorithm has a training phase that helps the customer to decide on the hyper-parameters that best describe their risk profile. The algorithm also has a validation (live production) phase that recommends the reorder parameters in real-time.

Note that the algorithm allows the customer to adjust the objective function based on their business situation. Also, the customers have the option to bypass the training phase and enter the live production with predefined hyper-parameters if they already have a good understanding of their risk profile. We pointed out that some idealized assumptions for the uncertainties in the literature oftentimes does not hold. The C3 SIO algorithm illustrated a practical workflow to handle the uncertainties in real-time production.

The performance of C3 SIO algorithm is demonstrated in detail with data provided by one of the clients in the energy and manufacturing industry. We also summarize the performance of the algorithm on customers' data across a variety of verticals. We showed that C3 SIO algorithm is effective in reducing the inventory and the safety stock costs while maintaining a satisfying service level.

A sensitivity study is conducted showing that some of the model parameters will impact the optimization results. The sensitivity study provided reasonable choices for model parameters as default values while deploying the algorithm in scale. Nevertheless, we still encourage the customers to test on a sub-sample of SKUs to find the best model parameters for their dataset before the live production.

Finally, we want to note here that C3 SIO algorithm is distinct from the solutions built with MILP in the literature due to its compatibility with the MRP systems, which carries critical value in terms of the customer adaptation. Should the MRP limitation be removed and the customers are open to adapt an end-to-end solution (i.e., instead of accepting reorder parameters, but rather the order placement recommendations), we would have more flexibility to utilize the classical solutions and further improve C3 SIO.

\section*{Acknowledgements}

We would like to thank C3.ai who provided all the resources that made this research possible. We would like to thank our customers to worked with us over the last few years, providing valuable feedback to help improving our C3 Inventory Optimization application. We would also like to thank our internal experts Alex Amato, Nikhil Krishnan, Burak Gundogdu for their valuable inputs. Our special thanks go to Stefano Zavagli who provided insightful comments and contributed to some of the figures in the paper.

\input{main.nls}

{\small
\bibliographystyle{apalike}

\begin{thebibliography}{}

\bibitem[Agarwal, 2019]{agarwal2019multi}
Agarwal, A. (2019).
\newblock Multi-echelon supply chain inventory planning using
  simulation-optimization with data resampling.
\newblock {\em arXiv preprint arXiv:1901.00090}.

\bibitem[Aikens, 1985]{AIKENS1985263}
Aikens, C. (1985).
\newblock Facility location models for distribution planning.
\newblock {\em European Journal of Operational Research}, 22(3):263--279.

\bibitem[Chu et~al., 2015]{chu2015simulation}
Chu, Y., You, F., Wassick, J.~M., and Agarwal, A. (2015).
\newblock Simulation-based optimization framework for multi-echelon inventory
  systems under uncertainty.
\newblock {\em Computers \& Chemical Engineering}, 73:1--16.

\bibitem[Fu, 2002]{10.5555/2700739.2700741}
Fu, M.~C. (2002).
\newblock Optimization for simulation: Theory vs. practice.
\newblock {\em INFORMS J. on Computing}, 14(3):192–215.

\bibitem[Garcia et~al., 1989]{garcia1989model}
Garcia, C.~E., Prett, D.~M., and Morari, M. (1989).
\newblock Model predictive control: Theory and practice—a survey.
\newblock {\em Automatica}, 25(3):335--348.

\bibitem[Geoffrion and Graves, 2010]{inbook}
Geoffrion, A. and Graves, G. (2010).
\newblock {\em Multicommodity Distribution System Design By Benders
  Decomposition}, volume~20, pages 35--61.

\bibitem[Geoffrion and Powers, 1995]{Geoffrion1995TwentyYO}
Geoffrion, A. and Powers, R.~F. (1995).
\newblock Twenty years of strategic distribution system design: An evolutionary
  perspective.
\newblock {\em Interfaces}, 25:105--127.

\bibitem[Gupta and Maranas, 2000]{Gupta2000ATM}
Gupta, A. and Maranas, C. (2000).
\newblock A two-stage modeling and solution framework for multisite midterm
  planning under demand uncertainty.
\newblock {\em Industrial \& Engineering Chemistry Research}, 39:3799--3813.

\bibitem[Ierapetritou and Pistikopoulos, 1994]{Ierapetritou1994NovelOA}
Ierapetritou, M. and Pistikopoulos, E. (1994).
\newblock Novel optimization approach of stochastic planning models.
\newblock {\em Industrial \& Engineering Chemistry Research}, 33:1930--1942.

\bibitem[{Joines} et~al., 2002]{1166395}
{Joines}, J.~A., {Gupta}, D., {Gokce}, M.~A., {King}, R.~E., and {Kay}, M.~G.
  (2002).
\newblock Supply chain multi-objective simulation optimization.
\newblock In {\em Proceedings of the Winter Simulation Conference}, volume~2,
  pages 1306--1314 vol.2.

\bibitem[Jung et~al., 2008]{JUNG20082570}
Jung, J.~Y., Blau, G., Pekny, J.~F., Reklaitis, G.~V., and Eversdyk, D. (2008).
\newblock Integrated safety stock management for multi-stage supply chains
  under production capacity constraints.
\newblock {\em Computers \& Chemical Engineering}, 32(11):2570--2581.
\newblock Enterprise-Wide Optimization.

\bibitem[Köchel and Nieländer, 2005]{KOCHEL2005505}
Köchel, P. and Nieländer, U. (2005).
\newblock Simulation-based optimisation of multi-echelon inventory systems.
\newblock {\em International Journal of Production Economics}, 93-94:505--513.
\newblock Proceedings of the Twelfth International Symposium on Inventories.

\bibitem[Mayne et~al., 2000]{mayne2000constrained}
Mayne, D.~Q., Rawlings, J.~B., Rao, C.~V., and Scokaert, P.~O. (2000).
\newblock Constrained model predictive control: Stability and optimality.
\newblock {\em Automatica}, 36(6):789--814.

\bibitem[Mele et~al., 2006]{mele2006simulation}
Mele, F.~D., Guillen, G., Espuna, A., and Puigjaner, L. (2006).
\newblock A simulation-based optimization framework for parameter optimization
  of supply-chain networks.
\newblock {\em Industrial \& Engineering Chemistry Research}, 45(9):3133--3148.

\bibitem[{Olafsson} and {Jumi Kim}, 2002]{1172871}
{Olafsson}, S. and {Jumi Kim} (2002).
\newblock Simulation optimization.
\newblock In {\em Proceedings of the Winter Simulation Conference}, volume~1,
  pages 79--84 vol.1.

\bibitem[Ptak and Smith, 2011]{ptak2011orlicky}
Ptak, C. and Smith, C. (2011).
\newblock {\em Orlicky's material requirements planning}.
\newblock McGraw-Hill Education.

\bibitem[Santoso et~al., 2005]{SANTOSO200596}
Santoso, T., Ahmed, S., Goetschalckx, M., and Shapiro, A. (2005).
\newblock A stochastic programming approach for supply chain network design
  under uncertainty.
\newblock {\em European Journal of Operational Research}, 167(1):96--115.

\bibitem[Silver, 1973]{silver1973heuristic}
Silver, E.~A. (1973).
\newblock A heuristic for selecting lot size quantities for the case of a
  deterministic time-varying demand rate and discrete opportunities for
  replenishment.
\newblock {\em Prod. Inventory Manage.}, 2:64--74.

\bibitem[Smith and Smith, 2014]{smith2014demand}
Smith, D.~A. and Smith, C. (2014).
\newblock {\em Demand driven performance: using smart metrics}.
\newblock McGraw-Hill Education.

\bibitem[Swaminathan et~al., 1998]{swaminathan1998modeling}
Swaminathan, J.~M., Smith, S.~F., and Sadeh, N.~M. (1998).
\newblock Modeling supply chain dynamics: A multiagent approach.
\newblock {\em Decision sciences}, 29(3):607--632.

\bibitem[Velasco~Acosta et~al., 2020]{velasco2020applicability}
Velasco~Acosta, A.~P., Mascle, C., and Baptiste, P. (2020).
\newblock Applicability of demand-driven mrp in a complex manufacturing
  environment.
\newblock {\em International Journal of Production Research},
  58(14):4233--4245.

\bibitem[Vidal and Goetschalckx, 1997]{VIDAL19971}
Vidal, C.~J. and Goetschalckx, M. (1997).
\newblock Strategic production-distribution models: A critical review with
  emphasis on global supply chain models.
\newblock {\em European Journal of Operational Research}, 98(1):1--18.

\bibitem[Wagner and Whitin, 1958]{wagner1958dynamic}
Wagner, H.~M. and Whitin, T.~M. (1958).
\newblock Dynamic version of the economic lot size model.
\newblock {\em Management science}, 5(1):89--96.

\bibitem[Zipkin, 2000]{zipkin2000foundations}
Zipkin, P.~H. (2000).
\newblock {\em Foundations of inventory management}.

\end{thebibliography}

}

\end{document}